\documentclass[11pt]{article}
\usepackage{geometry}                % See geometry.pdf to learn the layout options. There are lots.
\usepackage{graphicx}
\usepackage{amssymb}
\usepackage{amsmath}
\usepackage{epstopdf}
\usepackage[all]{xypic}

\newtheorem{lemma}{Lemma}[section] 
\newtheorem{propos}[lemma]{Proposition}
\newtheorem{example}[lemma]{Example}

\newtheorem{cor}[lemma]{Corollary}

\newtheorem{defin}[lemma]{Definition}
\newtheorem{remark}[lemma]{Remark}

\newcommand{\CF}{\hbox{{$\mathcal F$}}}

\newcommand{\CM}{\hbox{{$\mathcal M$}}}

\newcommand{\CO}{\hbox{{$\mathcal O$}}}
\newcommand{\cO}{\hbox{{$\mathcal O$}}}

\newcommand{\C}{\mathbb{C}}
\newcommand{\PP}{\mathbb{P}}
\newcommand{\R}{\mathbb{R}}
\newcommand{\RR}{\mathbb{R}}

\newcommand{\Z}{\mathbb{Z}}
\newcommand{\N}{\mathbb{N}}

\newcommand{\Spec}{\mathrm{Spec}}
\newcommand{\GKdim}{\mathrm{GK}}

\newcommand{\Hom}{\mathrm{Hom}}

\newcommand{\Ext}{\mathrm{Ext}}

\newcommand{\rank}{\,\mathrm{rank}\,}

\newcommand{\pd}{\partial}

\newcommand{\extd}{\mathrm{d}}

\newcommand{\tens}{\mathop{\otimes}}

\newcommand{\id}{\mathrm{id}}

\newcommand{\End}{\mathrm{ End}}

\newcommand{\Mod}{{\sf Mod}}

\newcommand{\pdol}{\overline{\partial}}
\def\D{{\Delta}}
\newcommand{\Projnc}{\mathrm{Proj_{nc}}}
\newcommand{\SL}{\mathrm{SL}}
\newcommand{\SU}{\mathrm{SU}}
\newcommand{\U}{\mathrm{U}}
\newcommand{\Qcoh}{{\sf  Qcoh}}
\newcommand{\Gr}{{\sf Gr }}

\newcommand{\Fdim}{{\sf Fdim }}
\newcommand{\QGr}{{\sf QGr }}

\newcommand{\G}{{\Gamma}}

\def\aol{{\bar a}}
\def\bol{{\bar b}}

\def\Mol{{\bar M}}

\def\cal{\mathcal}

\def\cL{{\cal L}}

\def\cO{{\cal O}}

\def\fM{{\mathfrak M}}

\def\fg{{\mathfrak g}}

\def\fl{{\mathfrak l}}
\def\fm{{\mathfrak m}}
\def\fn{{\mathfrak n}}
\def\fp{{\mathfrak p}}

\def\a{\alpha}
\def\b{\beta}

\def\l{\lambda}

\def\s{\sigma}

\def\ve{\varepsilon}

\def\lambdaol{{\bar \lambda}}

\def\xiol{{\bar \xi}}

\def\nablaol{{\overline \nabla}}
\def\th{\mathrm{th}}

\def\ve{{\varepsilon}}

\newenvironment{pf}{\noindent{\bf Proof.}}{\hfill $\square$\medskip}

\def\hdot{{\:\raisebox{2pt}{\text{\circle*{1.5}}}}}    % small bullet for complexes and grading

\title{Noncommutative complex differential geometry}
\author{Edwin Beggs and S.\ Paul Smith}
\begin{document}

\maketitle
%\section{}
%\subsection{}

\abstract{This paper defines and examines the basic properties of 
noncommutative analogues of almost complex structures, integrable 
almost complex structures, holomorphic curvature, cohomology, and holomorphic sheaves. The starting point is a differential structure on a noncommutative algebra defined in terms of a differential graded algebra.
This is compared to current ideas on noncommutative algebraic geometry.}

\section{Introduction}

\subsection{Philosophy}
This paper is about noncommutative complex analytic manifolds and holomorphic sheaf cohomology theory.  The classical theory of complex manifolds begins with 
a smooth manifold of even dimension
endowed with an `almost complex structure'. The Newlander-Nirenberg condition
\cite{NewNir} says when this almost complex structure actually comes from a complex coordinate system. 

Every smooth complex algebraic variety is a complex manifold. Kodaira's embedding theorem 
characterizes the compact complex manifolds that may be embedded in $\C\PP^n$. Chow's Theorem shows that every compact complex submanifold of $\C\PP^n$ is a smooth complex projective algebraic variety. 
Serre's GAGA, the abbreviation of the title \textit{
G\'eom\'etrie alg\'ebrique et g\'eom\'etrie analytique} of \cite{SerreGAGA}, shows that the algebraic and analytic properties of a smooth complex subvariety of $\C\PP^n$ are ``the same''. 

In sum, these results show that the basic definitions in complex differential and complex algebraic geometry are compatible. This harmony confirms the appropriateness of the basic definitions 
in the two fields. Hodge theory provides further compatibilities between complex differential and complex algebraic geometry. 

There is nothing like this in noncommutative geometry. There is much less certainty
about the `correct' definitions. This is entirely reasonable: there are more
noncommutative algebras than commutative ones, and some noncommutative algebras
can be really horrible. There are many ways in which ideas from commutative 
geometry might be carried over into the noncommutative world, and there have been many attempts to do this. But which ideas correspond to `geometry' rather than the general theory of noncommutative algebras? 

These ideas, whatever they are, must have a large and diverse collection of examples.
They should usually reduce to the corresponding classical structures, but we should not exclude 
ideas that only have non-trivial meaning in the noncommutative world. And most classical geometric ideas should have noncommutative analogues. Certainly,
in mathematical physics, it would make no sense if we were to say that just because the real world should be noncommutative because of quantum theory, we should no longer be allowed to teach general relativity because the concepts in it cannot really work. If the universe really is `noncommutative', then there probably are ideas in noncommutative geometry which reduce to, for example,
 geodesics and parallel transport in a classical limit. 

This, then, is the philosophy of the paper: There \textit{ought}  to be a concept of 
noncommutative complex analytic manifolds within which one should be able to carry
out many of the classical operations. To be meaningful, there  \textit{should} be a good number of examples. 
And, for the future, there should be the  \textit{possibility} of using this as a bridge between noncommutative complex analytic manifolds and noncommutative complex
algebraic varieties.

\subsection{Overview of the paper}
This paper deals with non-commutative analogues of almost complex structures, integrable 
almost complex structures, holomorphic curvature, cohomology, and holomorphic sheaves. 
We attempt to give minimal conditions for a noncommutative complex structure that
allow the cohomology of holomorphic sheaves to be constructed.  We will always refer to smooth or differentiable manifolds as {\it real} manifolds so as to distinguish them
from complex manifolds. 

An almost complex structure on a  non-commutative real manifold  
having a non-commutative $*$-algebra $A$ of ``$\C$-valued differentiable functions'' 
and a de Rham complex $(\Omega^\hdot A,\extd,*)$ 
with complex coefficients (a differential graded algebra)  will be defined as an $A$-module homomorphism $J:\Omega^1 A\to \Omega^1 A$ whose square is minus the identity, and satisfies
a reality condition. This data, together with an extension of $J$ to a derivation on 
$\Omega^\hdot A$, leads to an $A$-bimodule decomposition 
$\Omega^\hdot A = \oplus_{p,q} \Omega^{p,q}A$ into $J$-eigenspaces such that
 $(\Omega^{p,q}A)^*=\Omega^{q,p}A$ and
$$
\Omega^nA=\bigoplus_{p+q=n} \Omega^{p,q}A
$$
for all $n$, and to associated operators
$$
\pd:\Omega^{p,q}A \to \Omega^{p+1,q}A
\qquad \hbox{and} \qquad
\qquad \pdol:\Omega^{p,q}A \to \Omega^{p,q+1}A\ .
$$
Motivated by the Newlander-Nirenberg condition we say that the almost complex structure is integrable if $\extd \Omega^{1,0}A \subset \Omega^{2,0}A \oplus \Omega^{1,1}A$. If the 
almost complex structure is  integrable  we show that
$$
\extd=\pd+\pdol, \qquad \pd^2=0, \qquad \pd\pdol+\pdol\pd=0, \qquad \pdol^2=0.
$$

Suppose $(\Omega^\hdot A,\extd,*,J)$ is an integrable almost complex structure on 
%(the non-commutative geometric object corresponding to)
 $A$. 
%In Section \ref{sect.hol.mods} we define 
A holomorphic $A$-module is a pair $(E,\nablaol)$ where $E$ is a left $A$-module and $\nablaol:E \to  \Omega^{0,1}A \otimes_A E$ is a $\pdol$-connection/operator. If $\nablaol^2=0$ we call  $(E,\nablaol)$ a holomorphic left $A$-module. We define a category ${\sf Hol}(A)$ of such modules. Associated to 
each $(E,\nablaol) \in {\sf Hol}(A)$ is a complex from which we define cohomology groups $H^{{}^\bullet}(E,\nablaol)$. Elements of $E$ play the role of continuous sections of a sheaf on the 
underlying non-commutative geometric object and elements of $H^0(E,\nablaol)= \ker\nablaol$ play the role of holomorphic sections.

This category ${\sf Hol}(A)$  is abelian if $\Omega^{1,0}A$ is flat as a right $A$-module. 
If all $\Omega^{p,q}A$ are flat right $A$-modules every short exact sequence in ${\sf Hol}(A)$ yields a long exact sequence of the cohomology groups (Proposition \ref{prop.les.cohom}).

 \subsection{Relation to non-commutative complex projective algebraic geometry}
 \label{ssect.ncag}

Let $R$ be a finitely generated connected graded $\C$-algebra of Gelfand-Kirillov dimension 
$n+1$. We assume $R$ is  left and right noetherian, Artin-Schelter regular, and a domain. 
We think of $R$ as a homogeneous coordinate ring (hcr) of an irreducible non-commutative 
smooth complex projective variety,
$\Projnc R$, of dimension $n$, that is defined  implicitly by declaring that 
the category of ``quasi-coherent sheaves'' on it is
$$
\Qcoh (\Projnc R):= \QGr(R)= \frac{\Gr R}{\Fdim R}
$$
where $\Gr R$ is the category of $\Z$-graded $R$-modules and $\Fdim R$ is its full subcategory
of modules that are the sum of their finite dimensional submodules and $\QGr(R)$ is the quotient category.

$\Projnc R$ is a purely algebraic construct. 
In order to treat $\Projnc R$ as a non-commutative differential geometric object on which one can 
do calculus one first needs an underlying real structure on $\Projnc R$, 
presumably defined in terms of a $*$-algebra $A$ of Gelfand-Kirillov dimension $2n$ and a de Rham complex $(\Omega^\hdot A,\extd)$ which is also a $*$-algebra.

%Following \cite{ATV1}, \cite{ATV2}, and \cite{AZ}, 
%we write $\cO_X$ for the image of $R$ in $\Qcoh X$ and {\it define}
 %$H^q(X,-) := \Ext^p_{\Qcoh X}(\cO_X,-)$.
 
\subsubsection{}
 Stafford and Van den Bergh's survey article \cite{SV}  shows that many ideas, 
 tools, and results, of projective algebraic geometry extend to the non-commutative setting in 
 a seamless and satisfying way. 

However, there are no non-commutative analogues of Chow's Theorem, the Kodaira embedding Theorem, or Serre's GAGA principle. In the classical commutative setting those results allow the 
application of complex-analytic, Hodge-theoretic, and K\"ahler geometric, methods to complex projective varieties, and conversely the methods of algebraic geometry apply to appropriate
complex manifolds.

\subsubsection{}
Until this gap, or chasm, is bridged non-commutative geometry will be without a 
union of algebraic and analytic methods.  If the gap were bridged non-commutative 
geometry would be greatly enhanced. Such a development might also lead to a deep connection between non-commutative algebraic geometry and non-commutative geometry as conceived of by Connes, i.e., based on operator algebras. At present there are few links between the two subjects. 
The main contacts between the two schools include the following: the work of Connes and Dubois-Violette on the non-commutative 3-spheres related to the 4-dimensional Sklyanin algebra; the work of Polishchuk and Schwarz on holomorphic bundles on the non-commutative 2-torus; the work of several people on homogeneous spaces for quantum groups. 

\subsubsection{}
Part of the problem is that non-commutative complex projective algebraic geometry deals only 
with holomorphic aspects of geometry and the objects appearing in non-commutative complex projective algebraic geometry have no underlying real structure: although many of these objects behave like smooth complex projective varieties, with few exceptions,
there is no underlying smooth non-commutative real manifold on which a complex structure is then imposed.

Polishchuk and Schwarz's work on holomorphic structures on non-commutative tori \cite{PS1} is 
the great exception to the last statement and their work illustrates the advantages of being able to impose a complex structure on a non-commutative real manifold.

\subsubsection{Examples of complex structures on some non-commutative projective varieties}

We end this paper by showing how the ideas in this paper apply to 
some particular non-commutative projective algebraic varieties. 
The simplest of these examples, in \S\ref{ssect.CPn-theta},  is 
$\C\PP^n_\theta:=\Projnc R_\theta$ where $\theta$ is a
skew-symmetric $(n+1)\times(n+1)$ matrix over $\RR$, $R_\theta$ is the free algebra $\C\langle z_0,\ldots,z_n\rangle$ modulo relations $z_\mu z_\nu = \l^{\mu\nu} z_\nu z_\mu$ and $\l^{\mu\nu}=e^{{\rm i}\theta_{\mu\nu}}$.

More substantial examples, in \S\ref{ssect.q.flag.vars}, 
are the quantum group analogues of  irreducible generalized flag manifolds 
studied by Heckenberger and Kolb, \cite{HK06} and \cite{HK07}. 

Further specializing, the quantizations $\C\PP_q^n$  have been examined in detail by various subsets of F. D'Andrea, L. Dabrowski,  M.\ Khalkhali, G.\ Landi, A. Moatadelro, and W.D.\ van Suijlekom  (\cite{DADL}, \cite{DAL}, \cite{KhLavS}, \cite{KM}, \cite{KM2}).   The point of view of framed and associated bundles
is discussed in \cite{buncomproj}. 

The reader should be aware that the notation $\C\PP^n_q$ is used in different ways by different communities.  Although there is a feeling that the same object is being discussed the object 
called $\C\PP^n_q$ by one community belongs to a different category than the object called 
$\C\PP^n_q$ by the other community. 

Based on these examples we make some speculations in \S\ref{sect.nca+dg} about 
how one might impose integrable almost complex structures on the kinds of non-commutative projective varieties that appear in the ``Artin-Tate-Van den Bergh-Zhang'' version of non-commutative projective algebraic geometry.

 \subsection{Relation to other work}
Non-commutative differential geometry began with Connes's 1985 IHES paper of that name \cite{AC_NCDG}. 
That paper, the birth of quantum groups, and Woronowicz's 1987 and 1989 papers, \cite{W1} and \cite{W2}, in which he developed a $*$-differential calculus for quantum groups, were the beginning of 
a quarter-century development of non-commutative calculus. 

Non-commutative calculus is an algebraic creation  based on a non-commutative  algebra 
that plays the role of, and is viewed as, either the smooth or holomorphic functions on an imaginary non-commutative real manifold or complex manifold.  
Most of this work has focused on developing either an analogue of the 
classical de Rham complex of smooth forms or an analogue of the holomorphic de Rham complex. 

Less has been done to develop a non-commutative calculus modeled on that for a complex manifold, i.e., one that begins with the de Rham complex $(\Omega^\hdot A,\extd)$ of smooth forms, extends the scalars to $\C$, and then using a non-commutative analogue of a complex structure obtains a decomposition $\extd = \pd+\pdol$ and an associated decomposition 
$\Omega^\hdot_\C A = \oplus_{p,q} \Omega^{p,q}A$ into $(p,q)$-forms.  

A recent paper along these lines by Khalkali, Landi, and van Suijklekom \cite{KhLavS} begins with a short history of non-commutative complex geometry. That paper and others, for example, \cite{bdr}, \cite{CHZ}, \cite{DADL}, \cite{DAL}, \cite{DS}, \cite{KM}, \cite{K}, \cite{Maj}, \cite{PS1}, \cite{SV2}, 
examine a range of interesting examples. Most of those examples fit into our framework. 

Here we develop a fairly robust framework for a calculus of $(p,q)$-forms based on an almost 
complex structure $J$. Absent the $J$-operator, the definition of the $(p,q)$-forms in any 
particular example appears somewhat ad hoc although 
because the non-commutative examples that have been examined are
closely modeled on classical commutative examples the guidance
provided by the classical case has suggested what the $(p,q)$-forms should be.

  The completely different approach by twisting in \cite{bramaj}
 is an example of the current theory, as the classical integrability condition twists
 into the noncommutative one.

 \subsection{Acknowledgement}
 The authors would like to thank the Isaac Newton Institute, where this work began during the program Noncommutative Geometry, 24 July - 22 December 2006. The authors are also grateful to the  Milliman Fund  of the University of Washington's Mathematics Department for funding the first author's visit 
 there in September 2007. 
 Edwin Beggs would also like to thank the Royal Society for a grant to allow him to give a lecture on this material `Noncommutative sheaves and complex structures' at the International Conference on K-theory, C*-algebras and topology of Manifolds, Chern Institute of Mathematics, Tianjin, 1-5 June, 2009. 
 The second author is grateful for support from the NSF (grant DMS-0602347).   The authors would like to thank Stefan Kolb for providing notes relating to the papers
 \cite{HK04,HK06,HK07}, which were used in writing Section~\ref{ssect.q.flag.vars}.

\section{Almost complex structures on $*$-algebras}

In this paper we develop the rudiments of a non-commutative complex differential geometry based on the
notion of $*$-algebras that applies to {\it some} non-commutative algebraic varieties.
Several issues must be addressed in doing this. First, most non-commutative projective varieties are defined globally in terms of their homogeneous coordinate rings. They have no underlying topological space, and there is no patching together of affine pieces. Second, complex differential geometry is developed by imposing additional structure on an underlying real manifold but most non-commutative projective varieties do not seem to have an underlying non-commutative analogue of 
a real manifold. 
This forces us to develop a formalism based on algebraic data that takes the place of the underlying real manifold. That underlying algebraic data consists of an algebra with a $*$-structure.

\subsection{$*$-algebras}

\begin{defin}
A {\sf $*$-structure} on an associative $\C$-algebra $A$ is a map $A \to A$, $a \mapsto a^*$ such that $a^{**}=a$, $(ab)^*=b^*a^*$, and $(\lambda a + \mu b)^* =  \overline{\lambda} a^* + \overline{\mu} b^*$ for all $a,b \in A$ and all $\lambda, \mu \in \C$. We then call $A$ a {\sf $*$-algebra}. We call an element $a\in A$ {\sf self-adjoint} or {\sf hermitian} if $a^*=a$. A {\sf $*$-homomorphism} between $*$-algebras is a $\C$-algebra
homomorphism $f$ such that $f(a^*)=f(a)^*$ for all $a$. 
\end{defin}

Throughout this paper $A$ denotes an associative not-necessarily-commutative $\C$-algebra. 
If the $*$-structure is ignored, elements of $A$ should be thought of as playing the role of holomorphic 
functions on a ``non-commutative complex manifold'' or regular functions on a ``non-commutative 
quasi-affine variety''. With the $*$-structure, elements of $(A,*)$ should be thought of as $\C$-valued functions on a ``non-commutative real manifold'' or ``non-commutative real algebraic variety''.  Often $(A,*)$ plays the role of a dense subalgebra of a C${}^*$-algebra. 

For example, elements of the Hopf algebra $\cO_q(\SL(2,\C))$  can be thought of as regular functions on the non-commutative complex algebraic quantum group $\SL_q(2,\C)$ or, when
a suitable $*$-Hopf structure is considered, as $\C$-valued polynomial functions on its compact real
 form $\SU_q(2)$.

The matrix algebra $M_n(\C)$ is always given the $*$-structure that sends a matrix to its conjugate transpose.
An $n$-dimensional {\sf $*$-representation} of a $*$-algebra $A$ is a $*$-homomorphism $\varphi:A \to M_n(\C)$. 

The free algebra $F=\C\langle x_1,\ldots,x_n,x_1^*,\ldots,x_n^*\rangle$ on $2n$ variables 
is a $*$-algebra with $(x_i)^*=x_i^*$, and $\lambda^*=\lambdaol$ for $\lambda \in \C$. If $A$ is a finitely generated 
$*$-algebra there is a surjective $*$-homomorphism $F \to A$ for a suitable $n$. 

If $B$ is a commutative $\RR$-algebra $A:=\C \otimes_\RR B$ becomes a $*$-algebra by declaring that  $(\lambda \otimes b)^*:=\lambdaol \otimes b$  for $\lambda \in \C$ and $b \in B$. We recover $B$ as the subalgebra of self-adjoint elements. A $*$-homomorphism $f:A \to \C$ restricts to a homomorphism of $\RR$-algebras $B \to \RR$ and, conversely, every $\RR$-algebra homomorphism  $B \to \RR$ extends to a unique 
$*$-homomorphism $A \to \C$. Every  $*$-homomorphism is obtained by such an extension so extension and restriction are mutually inverse bijections
 between $*$-homomorphisms $A \to \C$ and  $\RR$-algebra homomorphisms  $B \to \RR$. 

These considerations apply when $B$ is the coordinate ring of a real algebraic variety. Suppose $X \subset \RR^n$ is the zero locus of a set of polynomials with real coefficients. Let 
$B=\RR[x_1,\ldots,x_n]/I$ where $I$
consists of the polynomials that vanish on $X$. Then $A:=\C \otimes B$ is a $*$-algebra and $X$ may be recovered as the set of $*$-homomorphisms $A \to \C$. We will write $\C[X]$ for 
$A$ with this $*$-structure. 
We may think of $\C[X]$ as the ring of $\C$-valued polynomial functions on $X$ with
$*$-structure defined by $f^*(x):=\overline{f(x)}$.

For example, if $S^1$ is the unit circle $x^2+y^2=1$ then $\C[S^1]$ is $\C[x,y]/(x^2+y^2-1)$ 
with $x^*=x$ and $y^*=y$. Thus, if $z=x+\mathrm{i}\,y$, then $z^*=x-\mathrm{i}\,y$ and $zz^*=1$ so 
$\C[S^1] \cong \C[z,z^*]/(zz^*-1)$. 
The unit circle is recovered as the $\R$-valued points for the subalgebra of self-adjoint elements
$\R[z+z^*,\mathrm{i}\,(z-z^*)] \cong \R[x,y]/(x^2+y^2-1)$. 

When $B$ is an $\RR$-algebra that is not commutative the identity map on $B$ is not a
$\C$-algebra anti-homomorphism so $\C \otimes_\RR B$ can't be given a $*$-structure by defining 
$(\lambda \otimes b)^*$ to be  $\lambdaol \otimes b$.  A standard example is 
provided by the algebra  $B=\R[x,x^*]$ generated by the annihilation and creation operators $x$ and $x^*$
 with relation $xx^*-x^*x=1$.  Now $A=\C \otimes_\RR B$ is the Weyl algebra.
 The momentum and  position operators $p$ and $q$ are the self-adjoint elements satisfying the equations
$$
x=\hbox{${{1}\over{2}}$}(q+\mathrm{i}p) \qquad \hbox{and} \qquad x^*=\hbox{${{1}\over{2}}$}(q-\mathrm{i}p).
$$
Then $pq-qp=-2\mathrm{i}$, illustrating the fact that the self-adjoint elements in a non-commutative $*$-algebra need not form a subalgebra.

\subsection{Conjugate bimodules}
\label{ssect.conj.bimods}
 
This mixing of differential geometry and $*$-structures must be done 
with a little care in noncommutative geometry - the reader is referred to \cite{barcats}
 for some of the details, most of which will not be necessary in this paper.

 Let $A$ be a $*$-algebra and $E$ an $A$-bimodule. We define the {\sf conjugate bimodule} $\overline{E}$ by declaring that 
 \begin{enumerate}
  \item 
   $\overline{E} =E$ as an  abelian group;
  \item 
  we write $\overline{e}$ for an element $e\in E$ when we consider it as an element of $\overline{E}$;
  \item 
   the bimodule operations for $\overline{E}$ are $a.\overline{e}=\overline{e.a^*}$
and $\overline{e}.a=\overline{a^*.e}$. 
\end{enumerate}
If $\theta:E \to F$ is any map we define   $\overline{\theta}:\overline{E}\to \overline{F}$ by
$\overline{\theta}(\overline{e}):= \overline{\theta(e)}$.  

If $\theta:E \to F$ is a homomorphism of $A$-bimodules, so is $\overline{\theta}$. 
In this way the bar operation $E\mapsto \overline{E}$ is a functor from the category ${}_A\mathcal{M}_A$ 
of $A$-bimodules to itself.

We make  $\overline{A}$ an associative algebra by defining the multiplication $\aol\bol:=\overline{ba}$. As an $\R$-algebra, $\overline{A}$ is isomorphic to $A^{\rm op}$ via the map $a \mapsto 
\overline{a}$.
We make  $\overline{A}$ a $\C$-algebra through the algebra homomorphism $\C \to \overline{A}$,
$\lambda \mapsto \overline{\lambda^*}$. We now define
\begin{equation}
\label{eq.star}
 \star:A\to \overline{A}, \qquad a\mapsto \overline{a^*}.
\end{equation}
Then $\star$  is a homomorphism of $\C$-algebras and an isomorphism  in ${}_A\mathcal{M}_A$.

\subsection{The universal differential calculus}

We adopt the standard notation and terminology for non-commutative differential calculus. The reader can find more details in \cite[Ch. 8]{GVF}, \cite[Ch. 12]{KS}, and \cite{AC_NCDG}.

Let $A$ be an arbitrary $\C$-algebra.  Let $\Omega^1_{\sf univ} A$ denote the kernel of the multiplication map $\mu:A \otimes A \to A$. Following Cuntz and Quillen \cite{CQ}, we define the universal differential graded algebra 
over $A$ to be the tensor algebra 
$$
\Omega^\hdot_{\sf univ} A = T_A(\Omega^1_{\sf univ} A)
$$ 
endowed with the unique degree one superderivation  such that
$$
\extd(a):=1 \otimes a - a \otimes 1
$$
for $a \in A$.

% We shall use symbols $\mathrm{i}$ for the square root of $-1$,  and $I$ for various identity operators.

\begin{defin}
A {\sf differential calculus} or {\sf differential structure} on $A$ is a differential graded algebra 
$(\Omega^\hdot A, \extd)$ that is a quotient of $(\Omega^\hdot_{\sf univ}A,\extd)$ by a differential 
graded ideal whose degree-zero component is zero. The cohomology of $(\Omega^\hdot A,\extd)$ is called the complex-valued {\sf de Rham cohomology} and denoted by $H^\hdot_{dR}(A)$.
\end{defin}

\noindent{\bf Warning:}
We will denote the product in $\Omega^\hdot A$ by wedge $\wedge$ although   
$\xi\wedge\eta$ will not usually equal  $(-1)^{|\xi|\,|\eta|}\eta\wedge\xi$ when $A$ is not commutative. 

\smallskip
It follows from the definition that 
   $\Omega^0A=A$; that $\Omega^{n+1}A=A.\extd(\Omega^n A)$ for all $n\ge 1$;
   that  the multiplication map
  $A \otimes \extd(\Omega^nA) \to \Omega^{n+1}A$ is surjective  for all $n \ge 0$; and
  $$\extd(\xi\wedge\eta)\,=\, \extd\xi\wedge\eta+(-1)^{|\xi|}\, \xi\wedge\extd\eta.$$
By hypothesis, $\Omega^\hdot A$ is generated by $A$ and $\Omega^1A$, and $\Omega^nA$ is the 
$\C$-span of 
 $$
\{ a_0\, \extd a_1 \wedge \cdots \wedge \extd a_n \; | \; a_0,\ldots,a_n \in A\}.
 $$
 
Let $(\Omega^\hdot A,\extd)$ and  $(\Omega^\hdot B,\extd)$ be differential calculi.
An algebra homomorphism $\phi:A \to B$ is {\sf differentiable} if it extends to a homomorphism
$\phi: \Omega^\hdot A \to \Omega^\hdot B$ of dgas. In particular, $\extd(\phi(a)) = \phi(\extd a)$ for all $a \in A$.

\subsection{Differential $*$-calculus}
 
Non-commutative $*$-calculus first appears in Woronowicz's papers 
\cite{W2} and \cite{W1} and \cite[Defn. 1.4]{W2} has been extended to higher degree forms as follows.

\begin{defin} \label{kjzxhdgv}
\cite[p.462]{KS}
A differential calculus $(\Omega^\hdot A,\extd)$ on a $*$-algebra $A$ is {\sf compatible with the star operation}
on $A$ if the star operation on $\Omega^0 A=A$ extends to an involution $\xi\mapsto\xi^*$ on $\Omega^\hdot A$ 
that preserves the grading and has the property that 
$(\extd\xi)^*=\extd(\xi^*)$ and $(\xi\wedge\eta)^*=(-1)^{\vert \eta \vert \vert \xi \vert}\eta^*\wedge\xi^*$ for all 
homogeneous $\eta, \xi\in\Omega A$. 

When these conditions hold we call $(\Omega^\hdot A,\extd,*)$  a {\sf differential $*$-calculus} on $A$. 
\end{defin}

\medskip 

\medskip\noindent{\bf Remark.}
The definition implies  $(a\xi b)^*=b^*\xi^* a ^*$ for all $a,b \in A$ and $\xi \in \Omega^\hdot A$. 

We can rephrase the definition of a differential $*$-calculus using the notion of a conjugate bimodule
in \S\ref{ssect.conj.bimods} and the fact that the product
\begin{eqnarray*}
\overline{\xi}\wedge\overline{\eta} &:=& (-1)^{|\xi|\,|\eta|}\,\overline{\eta\wedge\xi}\ 
\end{eqnarray*}
makes $(\overline{\Omega^\hdot A},\overline{\extd})$ a differential graded algebra.

\begin{propos}
\label{prop.*-calculus}
Suppose $(\Omega^\hdot A,\extd,*)$ is a differential $*$-calculus on $A$. 
\begin{enumerate}
  \item 
 If $\sum a_i.\extd b_i=0$ in $\Omega^1A$, then  $\sum_i \extd b_i^*.a_i^*=0$. 
 \item
 The map $\star:\Omega^\hdot A\to \overline{\Omega^\hdot A}$, $\star\xi:=\overline{\xi^*}$, is an $A$-bimodule homomorphism
 and a homomorphism of differential graded algebras.
 \end{enumerate}
 Conversely, if $\xi\mapsto\xi^+$ on $\Omega^\hdot A$ is an involution that extends the $*$ on $A$ and preserves the grading and has the property that $(\extd\xi)^+=\extd(\xi^+)$, then $(\Omega^\hdot A,\extd,{}^+)$ is a differential $*$-calculus on $A$. 
\end{propos}
\begin{pf}
We leave this to the reader. 
\end{pf}

 The map $\star:\Omega^\hdot A\to \overline{\Omega^\hdot A}$ is a morphism in the category of $A$-bimodules. The map $\Omega^\hdot A\to \Omega^\hdot A$, $\xi \mapsto \xi^*$ is not.

Our $\star$ has nothing to do with the Hodge star operation in Riemannian geometry, which in general
changes the degree of the form. Our $\star$ is a notational device that depends only on the fact that $A$ is a $*$-algebra. In 
the context of \textit{real} classical geometry $\star$ is the identity.

\begin{propos}
\label{prop.free.alg.*.calc}
Let  $F=\C\langle x_1,\ldots,x_n,x_1^*,\ldots,x_n^*\rangle$ be the free $*$-algebra. Then there is a differential 
$*$-calculus $(\Omega^\hdot_{\sf univ} F,\extd,*)$ defined by $(\extd a)^*=\extd(a^*)$ for $a \in F$ and
\begin{equation}
\label{eq.star.univ}
(a_0 \extd a_1 \cdots \extd a_n)^*:= \ve_n (\extd a_n^*) \cdots (\extd a_1^*)a_0^*
\end{equation}
where $\ve_n=(-1)^{{{1}\over{2}}n(n-1)}$; i.e., $\ve_n=1$ if $n$ is 0 or 1 modulo 4, and -1 otherwise. 
\end{propos}
\begin{pf}
 If $\xi=a_0\extd a_1 \ldots \extd a_n$, then
\begin{eqnarray*}
(\extd \xi)^* &=& (\extd a_0\extd a_1 \ldots \extd a_n)^*
\\
&= &\ve_{n+1} \extd a_n^* \cdots \extd a_1^* \extd a_0^*
\\
& =& \ve_n (-1)^n   \extd a_n^* \cdots \extd a_1^* \extd a_0^*
\\
&=& \ve_n \extd (\extd a_n^* \cdots \extd a_1^* a_0^*) 
\\
& =& \extd (\xi^*). 
\end{eqnarray*}
Every element in $\Omega^n_{\sf univ} F$ is a sum of terms of the form $a_0\extd a_1 \ldots \extd a_n$
so $(\extd \xi)^* = \extd (\xi^*)$ for all $\xi \in \Omega^n_{\sf univ} F$. 

If $a \in F$ and $\xi \in \Omega^\hdot_{\sf univ} F$ it follows at once that $(a\xi)^*=\xi^*a^*= a \overline{\xi^*}$
so $\star:\Omega^\hdot F\to \overline{\Omega^\hdot F}$ is a homomorphism of left $F$-modules. 

A straightforward calculation shows that
\begin{eqnarray*}
\extd a_0\extd a_1 \ldots \extd a_{n-1} a_n &=&  (-1)^n a_0 \extd a_1 \cdots \extd a_n 
\\
& & \qquad + \sum_{i=0}^{n-1}(-1)^{n-i+1} \extd a_0 \cdots \extd(a_ia_{i+1}) \cdots \extd a_n.
\end{eqnarray*}
It follows that $(\extd a_0\extd a_1 \ldots \extd a_{n-1} a_n)^*$ is equal to
\begin{eqnarray*}
 &&  (-1)^n \ve_n  \extd a_n^* \cdots \extd a_1^* a_0^*  + \sum_{i=0}^{n-1} (-1)^{n-i+1} \ve_n  \extd a_n^* \cdots \extd(a_ia_{i+1})^* \cdots \extd a_0^*.
\\
&=&  (-1)^n \ve_n \Big( (-1)^n a_n^* \extd a_{n-1}^* \cdots  \extd  a_0^* +  \sum_{i=0}^{n-1} (-1)^{n-i+1}   \extd a_n^* \cdots \extd(a_{n-i}^*a_{n-i-1}^*) \cdots \extd a_0^* \Big) 
\\
& \quad& 
+ \sum_{i=0}^{n-1} (-1)^{n-i+1} \ve_n  \extd a_n^* \cdots \extd(a_{i+1}^*a_{i}^*) \cdots \extd a_0^*
\end{eqnarray*}
The two  sums cancel so we obtain
$$
(\extd a_0\extd a_1 \ldots \extd a_{n-1} a_n)^* = a_n^* \extd a_{n-1}^* \cdots  \extd  a_0^*
$$
which shows that  $\star:\Omega^\hdot F\to \overline{\Omega^\hdot F}$ is a homomorphism of right $F$-modules. 
\end{pf}

We note that $\ve_n$ is the sign of the permutation $\Big( {1 \atop n}{2 \atop n-1} \cdots {n \atop 1}\Big)$
and (\ref{eq.star.univ}) therefore agrees with the formula in \cite[(3.3)]{W1}.  

\subsection{Almost complex structures}

An almost complex structure on a real manifold is an endomorphism  $J$ of its tangent bundle such that $J^2=-1$. A complex manifold, when viewed as a real manifold, has a natural almost complex structure \cite[Prop. 2.6.2]{Huy}. Since we are working with (non-commutative analogues of) differential forms rather than the tangent bundle we must express the notion of an almost complex structure in terms of differential forms.

\begin{defin}
\label{defn.acs}
Let $(\Omega^\hdot A,\extd,*)$ be a differential $*$-calculus on $A$. An {\sf almost complex structure} on 
 $(\Omega^\hdot A,\extd,*)$ is a degree zero derivation $J:\Omega^\hdot A \to \Omega^\hdot A$
 %(i.e.\ it obeys the same rules as $h$ in subsection \ref{cvshvcsv}) $J:\Omega^\hdot A \to \Omega^\hdot A$ 
 such that 
 \begin{enumerate}
 \item
 $J$ is identically 0 on $A$ and hence an $A$-bimodule endomorphism of $\Omega^\hdot A$;
  \item
  $J^2=-1$ on  $\Omega^1A$; and 
   \item 
$J(\xi^*)=(J\xi)^*$ for $\xi\in\Omega^1 A$ (or, equivalently,  $\overline{J}\,\star=\star\, J$ on $\Omega^1A$).
\end{enumerate}
\end{defin}

At this stage, there is no requirement on how $J$ interacts with $\extd$. Such a requirement will appear 
later when we define what it means for an almost complex structure to be integrable.

Because $J^2=-1$ on $\Omega^1A$, there is an $A$-bimodule decomposition
\begin{equation}
\label{10+01}
\Omega^1A = \Omega^{1,0}A \oplus \Omega^{0,1}A
\end{equation}
where  $\Omega^{1,0}A= \{ \omega \in \Omega^1A \; | \; J \omega ={\rm i}\,\omega\}$ 
and $\Omega^{0,1}A= \{ \omega \in \Omega^1A \; | \; J \omega =-{\rm i}\,\omega\}$.
Condition (3) in the definition implies $(\Omega^{0,1} A)^* = \Omega^{1,0} A$. 

\bigskip

The next result is analogous to the trivial fact that a differentiable manifold with an almost 
complex structure has even dimension, or, equivalently, the rank of its cotangent bundle
is even.

\begin{lemma}
\label{lem.rank}
Let $(\Omega^\hdot A,\extd,*)$ be a differential $*$-calculus on $A$. 
Suppose there is a function 
$$
\rank:\{\hbox{some $A$-bimodules}\} \to \Z
$$
defined on a class of $A$-bimodules that is closed under finite direct sums and direct summands and contains $A$ and all $\Omega^nA$ and has the property that   $\rank M = \rank \Mol$ and 
$\rank (M \oplus N)=\rank M + \rank N$.  If $A$ has an almost complex structure, then 
$\rank \Omega^1A$ is even.
\end{lemma}
\begin{pf}
The map $\star:\Omega^{1,0}A \to \overline{\Omega^{0,1}A}$ is an isomorphism of $A$-bimodules so
$\Omega^{1,0}A$ and $\Omega^{0,1}A$ have the same rank.  Therefore
$\rank \Omega^1A = 2 \rank \Omega^{0,1}A$. 
\end{pf}

\begin{remark}
The hypothesis  in Lemma \ref{lem.rank} is mild.  
Sometimes there will be a finite set that is a basis for $\Omega^1A$  as both a left and right $A$-module
and when $\Omega^1A$ is not free it is often a finitely generated projective $A$-module on the left and on the right, and becomes free after a suitable localization. 

It would be desirable to consider examples having the following additional properties: 
$A$ is a finitely presented $\C$-algebra, a domain, and noetherian, or coherent; $\Omega^1A$ is a finitely generated projective $A$-module of rank $2r$ on both the left and the right; %$\Omega^nA$ is a finitely generated projective $A$-module of rank ${2r}\choose{n}$ on both the left and the right; 
$\Omega^{p,q}A$, which we define in Lemma \ref{lem.tech} below, 
 is a finitely generated projective $A$-module of rank ${{r}\choose{p}}{{r}\choose{q}}$ on both the left and the right;
there is a left $A$-module isomorphism 
$\int_{\ell}:\Omega^{2r} A \to A$ and a right $A$-module isomorphism $\int_r:\Omega^{2r} A \to A$.
\end{remark}

\begin{propos}
\label{prop.free.alg.almost.cplx.str}
Let  $F=\C\langle x_1,\ldots,x_n,x_1^*,\ldots,x_n^*\rangle$ be the free $*$-algebra endowed with the  differential $*$-calculus $(\Omega^\hdot_{\sf univ} F,d,*)$ in Proposition \ref{prop.free.alg.*.calc}. Then 
there is a unique almost complex structure on $(\Omega^\hdot_{\sf univ} F,d,*)$ such that
$$
J(\extd x_k)=-\extd x_k^* \qquad \hbox{and} \qquad J(\extd x_k^*)=\extd x_k
$$
for $k=1,\ldots,n$
\end{propos}
\begin{pf}
Since $\Omega^1_{\sf univ} F$ is the free $F$-bimodule with basis 
$\{\extd x_k, \extd x_k^* \; | \; 1 \le k \le n\}$
any action of $J$ on the basis extends in a unique way to an $F$-bimodule automorphism of 
$\Omega^1_{\sf univ} F$. Then, because $\Omega^\hdot_{\sf univ} F$ is the tensor algebra over $F$ on 
$\Omega^1_{\sf univ} F$, the action of $J$ on $\Omega^1_{\sf univ} F$ extends in a unique way to an 
action of $J$ as a derivation of $\Omega^\hdot_{\sf univ} F$.  
\end{pf}

\noindent
{\bf Warning.}
When $n > 1$ the map $J:\Omega^nA \to \Omega^nA$ does not satisfy the equation 
$J^2=-1$. For example, because $J$ is a derivation $J(\extd a \wedge \extd b) = J(\extd a) \wedge \extd b+ \extd a \wedge J(\extd b)$ 
for all $a,b \in A$ and therefore  
\begin{eqnarray*}
J^2(\extd a \wedge \extd b) & = & J^2(\extd a)\wedge \extd b + 2 J(\extd a) \wedge J(\extd b)+\extd a \wedge J^2(\extd b) 
\\
& = &  2\big( J(\extd a) \wedge J(\extd b) - \extd a \wedge \extd b \big).
\end{eqnarray*}
More succinctly,
\begin{equation}
\label{J^2.Omega^2}
J^2=2(J \wedge J-1):\Omega^2A \to  \Omega^2A.
\end{equation}

\subsubsection{$(p,q)$-forms}
\label{sect.pq.forms}

Let $\C[h]$ denote the commutative polynomial ring viewed as the enveloping algebra of the 1-dimensional 
Lie algebra. Thus, if $V$ and $W$ are $\C[h]$-modules,  $V \otimes_\C W$ is given the $\C[h]$-module structure $h\cdot ( v \otimes w) =  v\otimes hw + hv\otimes w $.\footnote{Equivalently, if 
$\D:\C[h] \to   \C[h] \tens \C[h]$ is the  $\C$-algebra homomorphism defined by 
$
\D(h):=1 \otimes h + h \otimes 1,
$
then $V \otimes W$ is naturally a $ \C[h] \tens \C[h]$-module and is 
made into a $\C[h]$-module via $\D$. Alternatively, $V$ and $W$ can be viewed as representations of the trivial Lie algebra $\C h$ and $V^{\otimes r} \otimes W^{\otimes n}$ is then made into a representation of
$\C h$ in the standard way.} 
In particular, the tensor algebra $T(V)$ becomes a $\C[h]$-module and $h$ 
acts on $T(V)$ as a derivation. 

\begin{lemma}
\label{lem.tech}
Let $V$ be a $\C[h]$-module annihilated  by $h^2+1$. Let $n \ge 1$. 
Suppose   $\Omega^n$ is a $\C[h]$-module quotient of $V^{\otimes n} $. Then 
 there is a $\C[h]$-module decomposition
$$
\Omega^n = \bigoplus_{p+q=n} \Omega^{p,q}
$$
where
$$
\Omega^{p,q}=\{\xi \in \Omega^{\otimes n} \; | \; h.\xi=(p-q)\mathrm{i} \xi \}. 
$$
Furthermore, if $\Omega^\hdot = \bigoplus_{n=0}^\infty \Omega^{n}$ is a quotient algebra of $T(V)$ by a 
$\C[h]$-stable graded ideal, then the multiplication in  $\Omega^\hdot$ is such that $\Omega^{p,q} \otimes 
\Omega^{p',q'} \to \Omega^{p+p',q+q'}$. 
% In particular, $\C \oplus_{p \ge 1}\Omega^{p,0}$ is subalgebra. 
In particular, $\Omega^{0,\hdot}$ and $\Omega^{\hdot,0}$ are subalgebras of $\Omega^\hdot$. 
\end{lemma}
\begin{pf}
By hypothesis,  there is a $\C[h]$-module direct sum decomposition $V = V^{1,0} \oplus V^{0,1}$ where $h$ acts as multiplication by $\mathrm{i}$ on $V^{1,0}$ and as multiplication by $-\mathrm{i}$ on $V^{0,1}$. 

If $W$ and $W'$ are $\C[h]$-modules annihilated by $h-\lambda$ and $h-\lambda'$ respectively, then $W \otimes W'$
is annihilated by $h-(\lambda+\lambda')$. 
Applying this observation inductively to $(V^{1,0} \oplus V^{0,1})^{\otimes n}$ it follows that
$$
V^{\otimes n} = \bigoplus_{p+q=n} V^{p,q}
$$
where $V^{p,q}=\{\xi \in V^{\otimes n} \; | \; h.\xi=(p-q)\,\mathrm{i} \xi \}$. This proves the lemma for
 $\Omega^n=V^{\otimes n}$. Since $V^{\otimes n}$ is a semisimple $\C[h]$-module, so is every quotient of it, and the set of $h$-eigenvalues for a quotient of $V^{\otimes n}$ is a subset of the  set of $h$-eigenvalues for $V^{\otimes n}$.
 
This completes the proof of the first part of the lemma. The second part involving $\Omega^\hdot$ is equally
easy. 
 \end{pf}

\subsubsection{}
\label{ssect.Omega.pq}
 
Let $J$ be an almost complex structure on $(\Omega^\hdot A,\extd,*)$. 

By hypothesis the multiplication $\Omega^1 A \otimes_A \Omega^nA \to \Omega^{n+1}A$ is surjective for all $n\ge 1$. 
%This implies that the multiplication map is induced by a (non-unique) surjective vector-space homomorphism $\Omega^1A  \otimes_\C \Omega^nA \to \Omega^{n+1}A$. 
Therefore $\C \oplus \Omega^{\ge 1}A$ is a quotient of the tensor algebra $T_\C(\Omega^1A)$ where the tensor is taken over $\C$.

We may apply Lemma \ref{lem.tech} with $V=\Omega^1A$ and $h$  acting as $J$ does because
$J^2\xi =-\xi$ for all $\xi \in \Omega^1 A$. 
Because $T_\C(\Omega^1A)$ is generated  by $\Omega^1A$ as a $\C$-algebra there is a unique extension of $h$ to a derivation on $T_\C(\Omega^1A)$.
Since $J$ is a derivation of $\Omega^\hdot A$ the action of $J$ on 
$\Omega^nA$ is induced by the action of $h$ on $(\Omega^1A)^{\otimes n}$.

For all $p,q \ge 0$ we define
$$\Omega^{p,q}A : = \{\xi \in \Omega^{p+q} A \; | \; J\xi=(p-q)\mathrm{i}\,\xi\}.$$
Elements in $\Omega^{p,q}A$ are called {\sf $(p,q)$-forms}. 
Because $J$ is a derivation on $\Omega^\hdot A$ it follows from Lemma \ref{lem.tech} that
$$
\Omega^nA = \bigoplus_{p+q=n} \Omega^{p,q}A
$$
for all $n$. Because $J$ vanishes on $A$, it is a homomorphism of $A$-bimodules and therefore
each $ \Omega^{p,q}A$ is an $A$-bimodule.
By the last part of Lemma \ref{lem.tech},
\begin{equation} 
\label{eq.wedge.mult}
\Omega^{p,q} A \wedge \Omega^{p',q'} A \subset \Omega^{p+p',q+q'} A
\end{equation}
and $\Omega^{0,\hdot}A$ and $\Omega^{\hdot,0}A$ are subalgebras.

\subsubsection{The $\pd$ and $\pdol$ operators.}
\label{ssect.d.dbar}

For all pairs of non-negative integers $(p,q)$, let 
$$
\pi^{p,q}:\Omega^{p+q} A \to \Omega^{p,q}A
$$
be the projections associated to the direct sum decomposition $\Omega^nA= \oplus_{p+q=n} \Omega^{p,q}A$.

\begin{defin}
Let
$$
\partial :\Omega^{p,q}A \to \Omega^{p+1,q}A
$$
be the composition $\pi^{p+1,q}\extd: \Omega^{p,q} A \to  \Omega^{p+1,q} A$. 
Let
$$
\pdol :\Omega^{p,q}A \to \Omega^{p,q+1}A
$$
be the composition $\pi^{p,q+1}\extd: \Omega^{p,q} A \to  \Omega^{p,q+1}A$. 
\end{defin}

Because $\Omega^1A$ is generated as a left $A$-module by $\{\extd a \; | \; a \in A\}$,  
$\Omega^{1,0}A$ is generated as a left $A$-module by $\{\pd a \; | \; a \in A\}$; likewise, 
$\Omega^{0,1}A$ is generated as a left $A$-module by $\{\pdol a \; | \; a \in A\}$.

\begin{propos}
\label{prop.J}
Let $J$ be an almost complex structure on $(\Omega^\hdot A,\extd,*)$.
 Then for all $\xi\in\Omega^\hdot A$ and for all $p,q$:
 \begin{equation}
\label{prop.J.1}
J(\xi^*)=(J\xi)^*
\end{equation} 
 \begin{equation}
\label{prop.J.2}
(\Omega^{p,q}A)^*=\Omega^{q,p}A
\end{equation} 
\end{propos}
\begin{pf}
Certainly (\ref{prop.J.1}) holds when $\xi \in A\oplus \Omega^1A$. 
We now argue by induction on degree. Suppose (\ref{prop.J.1}) holds when $\xi \in \Omega^nA$. 
Every $n+1$ form is a sum of  elements $\xi\wedge\eta$ where
$\xi\in\Omega^nA$ and $\eta\in\Omega^1 A$. 
Because $J$ is a derivation on $\Omega^\hdot A$, 
\begin{eqnarray*}
(J(\xi\wedge\eta))^* &=& (J\xi \wedge\eta+\xi\wedge J\eta)^* \cr
&=& (-1)^n\, \big( \eta^*\wedge (J\xi)^* +  (J\eta)^*\wedge \xi^*\big) \cr
&=& (-1)^n\,\big( \eta^*\wedge J(\xi^*) +  J(\eta^*)\wedge \xi^*  \big)\cr
&=& (-1)^n\, J(\eta^*\wedge\xi^*) \cr
&=& J((\xi\wedge\eta)^*).
\end{eqnarray*}
Hence (\ref{prop.J.1}) holds when $\xi \in \Omega^{n+1}A$, and therefore holds for all $\xi \in \Omega^\hdot A$. Finally, a simple calculation shows that (\ref{prop.J.2}) follows from (\ref{prop.J.1}).
\end{pf}

\section{Integrable complex structures} 

\subsection{The Newlander-Nirenberg integrability condition}
By definition, an almost complex structure $J$ on a smooth manifold $X$ is a vector bundle isomorphism $J:TX \to TX$ such that $J^2=-1$. Such a $J$ leads to a decomposition $T^{1,0}X \oplus T^{0,1}X$ of the complexified tangent bundle $TX \otimes \C$ into the $\pm \mathrm{i}$-eigenspaces for $J$.
A complex manifold has a canonical almost complex structure and 
one says that $(X,J)$ is {\sf integrable} if there is a complex structure on $X$ such that $J$ is the canonical almost complex structure.  

 The Newlander-Nirenberg theorem says $J$ is integrable if and only if 
 $$
 [T^{0,1}X,T^{0,1}X] \subset T^{0,1}X.
 $$
 In terms of differential forms, the Newlander-Nirenberg theorem  says 
 $J$ is integrable if and only if 
 $$
 \extd\Omega^{1,0}_X \subset \Omega^{2,0}_X \oplus \Omega^{1,1}_X
 $$ 
 (see, e.g., \cite[Prop. 2.6.15]{Huy} or \cite[p. 54]{V}). Furthermore, 
$J$  is integrable if and only if $\extd \omega=\pd \omega+\pdol\omega$ for all $\omega \in \Omega^\hdot$  \cite[Prop. 2.6.15]{Huy}  if and only if $ \pdol^2=0$ on $A$  \cite[Cor. 2.6.18]{Huy}.

\begin{defin}\label{compdef}
An almost complex structure $J$ on $(\Omega^\hdot A,\extd,*)$ 
 is {\sf integrable} if any of the equivalent conditions in Lemma \ref{lem.int} hold (cf., \cite[Prop. 2.6.15 and Defn. 2.6.16]{Huy}).
\end{defin}

\begin{lemma}
\label{lem.int}
Let $J$ be an almost complex structure on $(\Omega^\hdot A,\extd,*)$.
The following conditions are equivalent:
\begin{enumerate}
\item{}
$\pdol^2 =0$ as an operator $A \to \Omega^2 A$;
\item{}
$\pd^2 =0$ as an operator $A \to \Omega^2 A$;
\item{}
$d=\pd +\pdol$ as operators $\Omega^1 A \to \Omega^2A$; 
  \item 
  $\extd \Omega^{1,0}A \subset \Omega^{2,0}A \oplus \Omega^{1,1}A$;
   \item 
  $\extd \Omega^{0,1}A \subset \Omega^{1,1} A\oplus \Omega^{0,2}A$.
\end{enumerate}
\end{lemma}
\begin{pf}
If $a \in A$, then $\extd^2 a=0$ so
$$
\pi^{2,0}\extd(\pd a + \pdol a) = \pi^{1,1}\extd(\pd a + \pdol a) = \pi^{0,2}\extd(\pd a + \pdol a) = 0.
$$
In other words,
\begin{equation}
\label{eq.zero.comps}
\pd^2 a + \pi^{2,0}\extd \pdol a = \pdol\pd a + \pd \pdol a = \pi^{0,2}\extd\pd a + \pdol^2 a = 0.
\end{equation}

(1) $ \Rightarrow$ (4)
Suppose (1) holds.  
Since $\pdol^2 a=0$, (\ref{eq.zero.comps}) implies $\pi^{0,2}\extd\pd a=0$, i.e., 
$\extd(\pd a) \in \Omega^{2,0}A \oplus \Omega^{1,1}A$. If $b \in A$, then 
$$
\extd(b. \pd a) = \extd b \wedge \pd a + b. \extd(\pd a) \in \Omega^1 A\wedge \Omega^{1,0} A \subset 
 \Omega^{2,0}A \oplus \Omega^{1,1}A\ .
 $$
 Since $\Omega^{1,0}A$ is generated as a left $A$-module by $\{b. \pd a \; | \; a,b \in A\}$ the last calculation shows that $\extd \Omega^{1,0}A \subset \Omega^{2,0}A \oplus \Omega^{1,1}A$.

(4) $\Rightarrow$ (1)
If $\extd \Omega^{1,0}A \subset \Omega^{2,0} A\oplus \Omega^{1,1}A$, then $\pi^{0,2} \,\extd(\pd a)=0$ for 
all $a \in A$ so $\pdol^2 a=0$ by  (\ref{eq.zero.comps}).

(2) $ \Leftrightarrow$ (5)
This is proved by the same kind of argument as was used to prove the equivalence of (1) and (4).

(4) $\Leftrightarrow$ (5)
Suppose (4) holds. Then $\big(\extd \Omega^{1,0} A\big)^* \subset \big(\Omega^{2,0}A\big)^* \oplus 
\big(\Omega^{1,1}A\big)^*$. Condition (5) now follows by applying  (\ref{prop.J.1}) and
  (\ref{prop.J.2}). 
The implication (5) $ \Rightarrow$ (4) is proved in a similar way.

(3) $\Rightarrow$ (4)
If $\extd=\pd+\pdol$ on $\Omega^1 A$, then 
$$
\extd \Omega^{1,0}A \subset \pd \Omega^1A+\pdol \Omega^{1,0} A\subset \Omega^{2,0} A+ \Omega^{1,1}A\ .
$$

(4) $ \Rightarrow$ (3)
Let $\omega \in \Omega^{1,0}A$ and $\eta \in \Omega^{0,1}A$. Then
$$
\extd \omega \in \extd \Omega^{1,0}A \subset   \Omega^{2,0} A+ \Omega^{1,1}A
$$
so $\extd \omega = \pd \omega+\pdol \omega$. Since (4) holds so does (5), and (5) implies that 
$\extd \eta \in \Omega^{1,1}A + \Omega^{0,2}A$ so  $\extd \eta = \pd \eta+\pdol \eta$. Thus $\extd(\omega+\eta)=(\pd+\pdol)(\omega+\eta)$. 
\end{pf}

\begin{lemma}
\label{lem.int2}
Let $J$ be almost complex structure on $(\Omega^\hdot A,\extd,*)$.
The following conditions are equivalent:
\begin{enumerate}
\item{}
$J$ is integrable;
  \item 
  $(1-J \wedge J)\extd J = J\extd:\Omega^1 A \to \Omega^2A$;
  \item
  $J^2\extd J = -2J \extd:\Omega^1A \to \Omega^2A$;
  \item 
   $J^2\extd  = 2 J \extd J :\Omega^1A \to \Omega^2A$;
   \item
   $J\extd J \extd =0:A \to \Omega^2A$.
\end{enumerate}
\end{lemma}
\begin{pf}
(1) $ \Rightarrow$ (3) and (4).
Because $\Omega^{1,0}A=(J+\mathrm{i})\Omega^1A$ and because $ \Omega^{2,0}A$, $\Omega^{1,1}A$, and  $\Omega^{0,2}A$ are the $2\mathrm{i}$-, $0$-, and $-2\mathrm{i}$-eigenspaces for $J$ acting on $ \Omega^{2}A$, Lemma \ref{lem.int}(4)
implies that $(J-2\mathrm{i})J\extd(J+\mathrm{i})=0$ on $\Omega^1A$. Lemma \ref{lem.int}(5) implies that
$(J+2\mathrm{i})J\extd(J-\mathrm{i})=0$ on $\Omega^1A$. But
$$
(J-2\mathrm{i})J\extd(J+\mathrm{i})= J^2\extd J+\mathrm{i}J^2 \extd -2\mathrm{i} J\extd J +2J\extd =0
$$
and 
$$
(J+2\mathrm{i})J\extd(J-\mathrm{i})= J^2\extd J- \mathrm{i}J^2 \extd +2\mathrm{i} J\extd J +2J\extd =0.
$$
It is now clear that (3) and (4) hold. 

(2) $ \Leftrightarrow$ (3)
We already showed (\ref{J^2.Omega^2}) that as endomorphisms of $\Omega^2A$,
$
J^2 = 2(J\wedge J -1).
$
Hence 
\begin{equation}
\label{eq.integ}
(1-J\wedge J) \extd J -J\extd = - \hbox{${{1}\over{2}}$} (J^2\extd J+2J\extd).
\end{equation}
The equivalence of (2) and (3) now follows. 

(3) $ \Leftrightarrow$ (4)
Since $J:\Omega^1 A \to \Omega^1 A$ is an invertible map, (3) holds if and only if (4) holds (multiply on the right by $J$). 

(4) $ \Rightarrow$ (1) 
Assume (4) holds. Then (3) holds too.

A 2-form $\xi$ belongs to $ \Omega^{2,0} A\oplus  \Omega^{1,1}A$ if and only if $(J+2\mathrm{i})J\xi=0$.
Since $\Omega^{1,0}A=(J+\mathrm{i})\Omega^1A$ it follows that Lemma \ref{lem.int}(4) holds if and only if $(J+2\mathrm{i})J\extd(J+\mathrm{i})$ vanishes identically on $\Omega^1A$. But 
$$
(J+2\mathrm{i})J\extd(J-\mathrm{i})= \big(J^2\extd J + 2J\extd \big) - \mathrm{i} \big( J^2 \extd - 2 J\extd J \big)
$$
and this is zero by (3) and (4).

(4) $ \Rightarrow$ (5) 
Since $J^2\extd = 2J\extd J$, $0=J^2\extd \extd= 2J\extd J \extd$. 

(5) $ \Rightarrow$ (4) 
Let $a,b \in A$. Then
$$
J^2\extd(a \extd b) = 2 (J \wedge J -1)(\extd a \wedge \extd b) = 2(J\extd a \wedge J\extd b -\extd a \wedge
 \extd b).
$$
On the other hand
$
2J\extd J(a\extd b)=2J(\extd a \wedge J\extd b +a \extd J \extd b) 
$
and $J(a\extd J \extd b)=aJ\extd J \extd b=0$, so 
$$
2J\extd J(a\extd b)=2J(\extd a \wedge J\extd b)=2(J\extd a \wedge J \extd b + \extd a \wedge J^2 \extd b).
$$
But $ \extd a \wedge J^2 \extd b =-\extd a \wedge \extd b$, so $J^2d=2JdJ:\Omega^1A \to \Omega^2A$.
\end{pf}

\begin{propos} \label{bvhdbvabv}
The operators 
$$
(1-J\wedge J)\extd J-J\extd, \quad J^2\extd J+2J\extd, \quad J^2\extd -2J \extd J
$$
are left $A$-module homomorphisms $\Omega^1 A \to \Omega^2A$. Hence $J$ is integrable if any one of these homomorphisms vanishes on a set of left $A$-module generators for $\Omega^1 A$.
(Similarly changing right for left everywhere.)
\end{propos}
\begin{pf}
Let us show that the first of these three operators is a left $A$-module homomorphism.
Let $a\xi \in \Omega^1 A$. Then 
\begin{eqnarray*}
\big( (1-J\wedge J)\extd J -J\extd\big)(a\xi) &=& (1-J\wedge J)\extd (a J \xi) -J(\extd a \wedge \xi + a \extd \xi) \cr
&=& (1-J\wedge J)(\extd a\wedge J\xi + a \extd J\xi)-J(\extd a \wedge \xi + a \extd \xi)\cr
&=& \extd a\wedge J\xi + J\extd a\wedge \xi+ a (1-J\wedge J)\extd J\xi -J(\extd a \wedge \xi + a \extd \xi) \cr
&=& a (1-J\wedge J)\extd J\xi  -  aJ\extd \xi.
\end{eqnarray*}
Hence $(1-J\wedge J)\extd J-J\extd$ is a left $A$-module homomorphism. It then follows from (\ref{eq.integ})
that the second operator is also a left $A$-module homomorphism. The third operator is obtained from the
second one by composing on the right with the $A$-module homomorphism $J$ so it too is a left $A$-module homomorphism.  
(The right module proof is similar.)
\end{pf}

\begin{propos}
\label{prop.pd+pdol}
Suppose  $(\Omega^\hdot A,\extd,*)$ is a differential $*$-calculus with an integrable almost complex structure  $J$. Then for all $p$ and $q$, 
$$
\extd \big(\Omega^{p,q}A\big)\; \subset \; \Omega^{p+1,q}A \oplus \Omega^{p,q+1}A.
$$
\end{propos}
\begin{pf}
We argue by induction on $p+q$. The case $p+q=1$ is true by the definition of 
integrability. Suppose  $n \ge 2$ and that the result is true for  all $p'+q'=n-1$. 

Let $p+q=n$. Then 
$$
\Omega^n A= \Omega^1 A\wedge \Omega^{n-1}A = \big(\Omega^{1,0}A\oplus \Omega^{0,1}A\big) 
 \wedge \Omega^{n-1}A
$$  
It therefore follows from (\ref{eq.wedge.mult}) that 
$$
\Omega^{p,q}A = \Omega^{1,0}A \wedge \Omega^{p-1,q}A + \Omega^{0,1}A \wedge \Omega^{p,q-1}A\ .
$$
However,  $\extd \Omega^{1,0}A \subset \Omega^{2,0} A\oplus \Omega^{1,1}A$ and 
  $\extd \Omega^{0,1} A \subset \Omega^{1,1} A\oplus \Omega^{0,2}A$ so, in conjunction with the induction 
  hypothesis applied to $ \extd  \Omega^{p-1,q}A$ and $ \extd  \Omega^{p,q-1}A$, we obtain
\begin{eqnarray*}
\extd \Omega^{p,q} &=&
 \extd \Omega^{1,0}A \wedge \Omega^{p-1,q} A\, + \, \Omega^{1,0}A \wedge  \extd  \Omega^{p-1,q}A
  \, + \, 
 \extd  \Omega^{0,1} A\wedge  \Omega^{p,q-1} A \, + \,   \Omega^{0,1} A \wedge  \extd   \Omega^{p,q-1} A
 \cr 
 &\subset &  \Omega^{p+1,q} A +  \Omega^{p,q+1} A \ ,
 \end{eqnarray*}
   as claimed.
\end{pf}

\begin{propos}
Suppose  $(\Omega^\hdot A,\extd,*)$ is a differential $*$-calculus with an integrable almost complex structure  $J$.  Then $\extd=\partial+\bar \partial$ and
\begin{equation}
\label{ddbar}
\partial^2=0, \qquad  \qquad \partial\bar \partial + \bar \partial \partial=0, \qquad \bar \partial^2=0.
\end{equation}
\end{propos}
\begin{pf}
By Proposition \ref{prop.pd+pdol}, $\extd=\partial+\bar \partial$. 
We have 
\begin{eqnarray*}
0 & = & \extd^2\Omega^{p,q}A
\cr
& \subset &  \,\partial^2 \Omega ^{p,q} A+ \big( \pd\pdol + \pdol\pd \big) \Omega ^{p,q} A
+ \pdol^2 \Omega ^{p,q} A
\cr
& \subset & \Omega ^{p+2,q}A \oplus \Omega ^{p+1,q+1}A \oplus \Omega ^{p,q+2}A\ . 
\end{eqnarray*}
The equalities in (\ref{ddbar})   follow  from the direct sum decomposition of $\Omega^nA$.
\end{pf}

\begin{propos}
Suppose  $(\Omega^\hdot A,\extd,*)$ is a differential $*$-calculus with an integrable almost complex structure  $J$. 
Both $\bar\partial$ and $\partial$ are super-derivations.
%, i.e.
%\begin{eqnarray*}
%\bar\partial(\xi\wedge\eta) &=& \bar\partial(\xi)\wedge\eta + (-1)^{|\xi|}
%\xi\wedge\bar\partial(\eta)\ ,\cr
%\partial(\xi\wedge\eta) &=& \partial(\xi)\wedge\eta + (-1)^{|\xi|}
%\xi\wedge\partial(\eta)\ .
%\end{eqnarray*}
\end{propos}
\begin{pf}
Let $\pi^{r,s}:\Omega^{r+s}A \to \Omega^{r,s}A$ denote the obvious projection (and its restrictions) for the 
direct sum decomposition of $\Omega^{r+s}A$ onto the direct sum of its subspaces of $(p,q)$-forms.

Let $\xi \in \Omega^{p,q}A$ and $\eta \in \Omega^{p',q'}A$. Then
\begin{eqnarray*}
\pd(\xi \wedge \eta) & = & \pi^{p+p'+1,q+q'}\extd(\xi \wedge \eta) 
\cr 
& = & \pi^{p+p'+1,q+q'}\big( \extd \xi \wedge \eta + (-1)^{|\xi|} \xi \wedge\extd  \eta \big) 
\cr 
& = & \big( \pi^{p+1,q} \extd \xi \big) \wedge \eta + (-1)^{|\xi|} \xi \wedge  \big( \pi^{p'+1,q'} \extd  \eta \big) 
\cr 
& = & \pd \xi  \wedge \eta + (-1)^{|\xi|} \xi \wedge  \pd  \eta
\end{eqnarray*}
as required. The proof for $\pdol$ is  essentially the same.
\end{pf}

 \begin{propos}
For all $\xi \in \Omega^*A$, 
$$
\pdol(\xi)^*=\pd(\xi^*)
\qquad \hbox{and} \qquad 
\pd(\xi)^*=\pdol(\xi^*).
$$
\end{propos}
\begin{pf}
Since it suffices to prove the proposition when $\xi \in \Omega^{p,q}A$ we make that assumption.
Since $\extd(\xi^*)=(\extd \xi)^*$, 
$$
(\pd \xi)^*+(\pdol\xi)^*=\pd(\xi^*) + \pdol(\xi^*).
$$
But $\xi^* \in \Omega^{q,p}A$ by Proposition \ref{prop.J}, 
so $=\pd(\xi^*) \in  \Omega^{q+1,p}A$ and $\pdol(\xi^*) \in  \Omega^{q,p+1}A$. 
The result now follows from the fact that $(\pd \xi)^* \in \big(\Omega^{p+1,q}A\big)^* =  \Omega^{q,p+1}A $ and 
 $(\pdol \xi)^* \in \big(\Omega^{p,q+1}A\big)^* =  \Omega^{q+1,p} A$. 
\end{pf}

\subsection{Holomorphic forms and holomorphic elements of $A$}
\label{sect.pd.pdol.1}

Let $J$ be an integrable almost complex structure on a differential
$*$-calculus $(\Omega^{\hdot}A,\extd,*)$.

An element $f \in A$ is {\sf holomorphic} if $\pdol f=0$.  
We define
\begin{equation}
A_{\sf hol}  := \{f \in A \; | \; \pdol f=0\}
\end{equation}
and
\begin{equation}
\label{defn.Omega.hol}
\Omega^p_{\sf hol} A :=\{ \omega \in \Omega^{p,0}A \; | \; \pdol \omega=0\}.
\end{equation}
Elements in $A_{\sf hol}$ are called {\sf holomorphic} (functions)
and elements in $\Omega^p_{\sf hol}A$ are called {\sf holomorphic p-forms}.

Elements in $A_{\sf hol}$ play the role of  ``holomorphic functions on a non-commutative complex variety'' and elements in $A$ play the role of ``$\C$-valued differentiable functions on an underlying 
real variety''.

As some reassurance, we check the analogue of the fact that the only $\RR$-valued holomorphic
 functions on a connected complex manifold are the constants. The analogue of an 
$\RR$-valued function is a self-adjoint element of $A$. 
 Suppose then that $f=f^*$ and $\pdol f=0$. Then 
$\extd f \in \Omega^{1,0}A$ so $(\extd f)^* \in \Omega^{0,1}A$. However, $(\extd f)^*=\extd(f^*)=\extd f \in \Omega^{1,0}A$ so  $\extd f \in  \Omega^{0,1}A \cap \Omega^{1,0}A$. Hence $\extd f = 0$. 
The analogue of connectedness
for a differential calculus is that
$\extd:A \to \Omega^1A$ vanishes only on $\C$; we therefore deduce that $f \in \C$, i.e., $f$ is constant.

\begin{propos}
\label{prop.Omega.hol.is.A.hol.bimod}
$A_{\sf hol}$ is a $\C$-subalgebra of $A$ and every $\Omega^p_{\sf hol}A$ is a bimodule over
$A_{\sf hol}$.
\end{propos}
\begin{pf}
Since $\pdol$ is a $\C$-linear derivation, $A_{\sf hol}$ is a $\C$-subalgebra of $A$.
If $f \in A_{\sf hol}$ and $\omega \in \Omega^p_{\sf hol}A$, then $\pdol(f\omega)=\pdol f.\omega + f \pdol \omega =0$ so $f \omega \in \Omega^p_{\sf hol}A$. Similarly, $\omega f \in \Omega^p_{\sf hol}A$.
\end{pf}

The  {\sf holomorphic de Rham complex} for $A$ is the complex
$$
0 \longrightarrow A_{\sf hol} \stackrel{\pd}{\longrightarrow} \Omega^1_{\sf hol} A\stackrel{\pd}{\longrightarrow} \Omega^2_{\sf hol} A \stackrel{\pd}{\longrightarrow} \cdots.
$$ 

\subsubsection{}
Omitting some crucial definitions, consider
a non-commutative projective algebraic variety $X:=\Projnc R$. For $\CF \in \Qcoh X$, 
one defines $H^q(X,\CF):=\Ext^q(\CO_X,\CF)$.
There is no general definition of objects $\Omega^p_X $ in $\Qcoh X$ though in some cases where 
$R$ is a Koszul algebra with properties like a polynomial ring there are reasonable candidates for
such $\Omega^p_X$ defined in terms of the Koszul resolution of the trivial $R$-module. It would be very interesting to examine whether there are situations in which $H_{\pdol}^{p,q}(A)$  is isomorphic to 
$H^q(X,\Omega_X^p)$ or the $q^{\th}$ cohomology group of the holomorphic de Rham complex for an appropriate $A$ endowed with a suitable integrable almost complex structure on $(\Omega^\hdot A,\extd,*)$.
 %and elements in $\Omega^{0,1}$ are 
%called {\sf anti-holomorphic 1-forms}. On $A$, $\pd:A \to \Omega^{1,0}$ and $\pdol:A \to \Omega^{0,1}$
%are given by
%\begin{equation}
%\label{eq.pd.pdol}
%\pd := \hbox{${{1}\over{2\mathrm{i}}}$}(J+\mathrm{i})d   \qquad \hbox{and } \qquad \pdol:= - \hbox{${{1}\over{2\mathrm{i}}}$}(J-\mathrm{i})d .
%\end{equation}
%Since $d:A \to \Omega^1A$ is a derivation and $J$ acts on $\Omega^1A$ as an $A$-module 
%homomorphism, $\pd$ and $\pdol$ are derivations. 

\section{Holomorphic modules}
\label{sect.hol.mods}

\subsection{The Koszul-Malgrange theorem}
Throughout this section $J$ is an integrable almost complex structure on 
 a differential $*$-calculus $(\Omega^\hdot A,\extd, *)$. In this section we will write  $\Omega^{p,q}A$ for the $A$-bimodules $\Omega^{p,q}A$.
 
Holomorphic $A$-modules are defined in Definition \ref{defn.holom.str} below. 
The justification for our definition is the  Koszul-Malgrange Theorem as we now explain. 

Let $E$ be a complex vector bundle on a complex manifold $M$. A holomorphic structure on $E$ is a complex manifold structure on the total space of $E$ such that the transition functions are 
holomorphic. Although there is no natural exterior derivative on a general complex vector bundle a {\it holomorphic} bundle $E$ does have a naturally defined $\C$-linear $\pdol$-operator 
$$
\pdol: \Omega^{p,q}A  \otimes_A E \to  \Omega^{p,q+1}A  \otimes_A E
$$
satisfying $\pdol^2=0$ and the Leibniz rule $\pdol(f\cdot \alpha) = \pdol(f) \wedge \alpha + f \pdol(\alpha)$ (see, e.g., \cite[p.70]{GH} and \cite[Lem. 2.6.23]{Huy}). 

Let $E$ be a complex vector bundle with a connection $\nabla:E \to \Omega^1A \otimes_A E$. 
We can write $\nabla=\nabla^{1,0}+\nabla^{0,1}$ where $\nabla^{p,q}:E \to \Omega^{p,q}A \otimes_A E$.
The Koszul-Malgrange Theorem, which is analogous to the Newlander-Nirenberg criterion, 
says that if the composition
$$
  \UseComputerModernTips
  \xymatrix{ 
  E \ar[rr]^{\nabla^{0,1}} && \Omega^{0,1} A \otimes_A E \ar[rr]^{\nabla^{0,2}} && \Omega^{0,2}A   \otimes_A E
  }
  $$
 is zero, in which case $\nabla^{0,1}$ is said to be integrable, then there is a unique holomorphic structure on $E$ with respect to which the natural $\pdol$-operator is $\nabla^{0,1}$. 

A given complex vector bundle $E$ might have many holomorphic structures each of which determines a
natural $\pdol$-operator which, in turn,  determines the holomorphic structure on $E$. 
Moreover, every $\C$-linear operator $\nablaol:E \to \Omega^{1,0} A \otimes_A E$ satisfying 
the conditions in Definition \ref{defn.holom.str} is the $\pdol$-operator for a unique holomorphic structure on  $E$ (see, e.g., \cite[Thm. 2.6.26]{Huy}).

\subsection{$\pdol$-operators on modules}

\begin{defin}
Let  $(\Omega^\hdot A,\extd,*)$ be a differential $*$-calculus with an integrable almost complex structure  $J$.
A {\sf $\pdol$-operator} on a left $A$ module $E$ is a $\C$-linear map
$$
\nablaol:E\to\Omega^{0,1} A \otimes_A E
$$
such that 
\begin{eqnarray*}
\nablaol(a.e) &=& \pdol a\otimes e + a.\nablaol e
\end{eqnarray*}
for all $a \in A$ and $e \in E$. For $q \ge 1$ we define
$$
\nablaol :\Omega^{0,q} A \otimes_A E\to \Omega^{0,q+1} A \otimes_A E
$$
 by
 \begin{equation}
\label{eq.hol.curv}
\nablaol(\xi\tens e)  \,:=\, \bar\partial \xi\tens e + (-1)^q\, \xi\wedge \nablaol (e)\ .
\end{equation}
The {\sf holomorphic curvature} of $\nablaol$ is defined to be
$$
\nablaol^2 :E\to \Omega^{0,2}A \tens E.
$$
\end{defin}

It is easy to check that 
\begin{enumerate}
  \item $\nablaol: \Omega^{0,q} A \otimes_A E\to \Omega^{0,q+1} A \otimes_A E$ is a homomorphism of left 
  $A_{\sf hol}$-modules and 
  \item 
$\nablaol^2$ is a left $A$-module map.
\end{enumerate}
% the holomorphic curvature  $R^{\nablaol}_{\mathrm{hol}}$

\begin{lemma}
\label{lem.nabla}
%Let  $(\Omega^\hdot A,\extd,*)$ be a differential $*$-calculus with an integrable almost complex structure  $J$.
Let $(E,\nablaol)$ be a module with a $\pdol$-operator. Then 
$$\nablaol^2:\Omega^{0,*} A\otimes_A E \longrightarrow   \Omega^{0,*+2}A \otimes_A E$$ is a homomorphism of left $\Omega^{0,*}A$-modules. In fact, for $\xi \in \Omega^{0,*}A$
and $e\in E$,
\begin{eqnarray*}
\nablaol^2(\xi \otimes e) &=& \xi \wedge      \nablaol^2(e).          %R^{\nablaol}_{\mathrm{hol}}(e)\ .
\end{eqnarray*}
\end{lemma}
\begin{pf}
An easy computation shows that $\nablaol^2(\xi \otimes e)=\xi \wedge \nablaol^2 e$.
Hence for every $\eta \in \Omega^{0,*}A$, 
$$
\nablaol^2(\eta \wedge \xi \otimes e)=\eta \wedge \xi \wedge \nablaol^2 e = \eta \wedge \nablaol^2(\xi \otimes e).
$$
In other words,  $\nablaol^2$ commutes with ``left multiplication by $\eta$'' so  is a homomorphism 
of left $\Omega^{0,*}A$-modules.
\end{pf}

\subsection{Holomorphic modules}
\label{ssect.hol.mods}

\begin{defin}
\label{defn.holom.str}
Let  $(\Omega^\hdot A,\extd,*)$ be a differential $*$-calculus with an integrable almost complex structure  $J$.
A {\sf holomorphic structure} on a left $A$-module $E$ is a
$\pdol$-operator $\nablaol$ on $E$ whose holomorphic curvature vanishes.   
%$R^{\nablaol}_{\mathrm{hol}}=0$. 
We then call $(E,\nablaol)$ a {\sf holomorphic}
$A$-module and an element $e \in E$ such that $\nablaol e=0$ is said to be  {\sf holomorphic}.
\end{defin}

\begin{defin}
The holomorphic modules $(E,\nablaol)$ are the objects in a category ${\sf Hol}(A)$. 
A morphism  $\phi:(E_1,\nablaol_1)\to (E_2,\nablaol_2)$ in ${\sf Hol}(A)$ is an $A$-module homomorphism  $\phi:E_1 \to E_2$ such that  $\nablaol_2 \phi = (\id \otimes \phi)\nablaol_1$.  
\end{defin}

\begin{propos}
If  $\Omega^{0,1} A$ is a flat right $A$-module, then ${\sf Hol}(A)$ is an abelian category.
\end{propos}
\begin{pf}
Let $f:(E_1,\nablaol_1) \to (E_2,\nablaol_2)$ be a morphism in ${\sf Hol}(A)$. Let $K$ be the kernel
and $C$ the cokernel of $f$ in the category of right $A$-modules. The hypothesis  that $\Omega^{0,1} A$ is flat implies the second row in the diagram
$$
 \UseComputerModernTips
  \xymatrix{ 
  0 \ar[r] & K \ar[r] &  E_1  \ar[d]_{\nablaol_1}  \ar[r]^f & E_2 \ar[r]  \ar[d]^{\nablaol_2}& C \ar[r] & 0
  \\
   0 \ar[r] &  \Omega^{0,1} A \otimes_A K \ar[r] &  \Omega^{0,1} A \otimes_A E_1   \ar[r]_{1 \otimes f} &  \Omega^{0,1} A \otimes_A E_2 \ar[r]  &  \Omega^{0,1} A \otimes_A C \ar[r] & 0
}
$$
is exact.
There are unique maps $\nablaol_K:K \to  \Omega^{0,1} A \otimes_A K$ and  $\nablaol_C:C \to  \Omega^{0,1} A \otimes_A C$ making the diagram commute. Since $\nablaol_K$ is, in effect, the restriction of $\nablaol_1$ it is a $\pdol$-operator. Likewise, $\nablaol_C$ is a $\pdol$-operator because it is induced by the $\pdol$-operator $\nablaol_2$. 

It is easy to check that  $(K,\nablaol_K)$ and $(C,\nablaol_C)$ are a kernel and 
cokernel for $f$ in ${\sf Hol}(A)$. 
\end{pf}

\medskip
If $(E,\nablaol)$ is a holomorphic module Lemma \ref{lem.nabla} implies that 
\begin{eqnarray*}
0 \longrightarrow E \stackrel{\nablaol} \longrightarrow \Omega^{0,1}A\otimes_A E 
 \stackrel{\nablaol} \longrightarrow \Omega^{0,2}A\otimes_A E 
  \stackrel{\nablaol} \longrightarrow \cdots
\end{eqnarray*}
is a complex. The cohomology groups of the complex  $\Omega^{0,\hdot} A\otimes_A E$ will be denoted by
$$
H^\hdot(E,\nablaol).
$$
If $\phi:(E_1,\nablaol_1)\to (E_2,\nablaol_2)$ is a homomorphism between holomorphic modules, then
$\id \tens \phi$ is a map of cochain complexes so  induces a map in cohomology 
$$
\phi^\hdot: H^\hdot(E_1,\nablaol_1) \to H^\hdot(E_2,\nablaol_2).
$$

 \begin{propos}
 \label{prop.les.cohom}
 Suppose every $\Omega^{n} A$ is flat as a right $A$-module. Then every exact sequence
 $$
 \UseComputerModernTips
  \xymatrix{ 
 0 \ar[r] &  (E,\nablaol_E) \ar[r]^{\phi} & (F,\nablaol_F)  \ar[r]^{\psi}  & (G,\nablaol_G) \ar[r] &  0
 }
 $$
 of holomorphic modules gives rise to a long exact sequence 
$$
\UseComputerModernTips
  \xymatrix{ 
 0 \ar[r] & H^0(E,\nablaol_E) \ar[r]^{\phi^*} & H^0(F,\nablaol_F) \ar[r]^{\psi^*} & H^0(G,\nablaol_G) \ar[r]^>>>>>>\delta & 
 \\
 & H^1(E,\nablaol_E) \ar[r]^{\phi^*} & H^1(F,\nablaol_F)  \ar[r]^{\psi^*} & \quad \cdots \quad
  }
 $$
of left $A_{\sf hol}$-modules.
 \end{propos}
 \begin{pf}
 A direct summand of a flat module is flat so each $\Omega^{0,n}A$ is flat as a right $A$-module.
 It follows that 
$$
\UseComputerModernTips
  \xymatrix{ 
 0 \ar[r] & \Omega^{0,\hdot}A \otimes_A E \ar[r] & \Omega^{0,\hdot} A \otimes_A F \ar[r] & \Omega^{0,\hdot} A \otimes_A G
 \ar[r] & 0
 }
 $$
 is an exact sequence of cochain complexes. 
 The existence of the long exact sequence follows in the usual way (see, e.g., \cite[Thm. 1.3.1]{Weib}).
 \end{pf}

 {\bf A Question.}
 Let $R$ be one of the non-commutative homogeneous coordinate rings that appears in non-commutative
 projective algebraic geometry over $\C$. Let $X = \Projnc R$. Let $M$ be a finitely generated graded left
  $R$-module that corresponds to an $\CO_X$-module $\CM$. Is there a relation between $\Ext^\hdot(\CO_X,\CM)$
  and $H^\hdot(E,\nablaol)$ for some holomorphic module $(E,\nablaol)$ over some $*$-calculus 
  $(\Omega^\hdot A,\extd,*)$ having an integrable structure  $J$?

\subsection{A connection whose $(0,2)$-curvature component vanishes induces a holomorphic structure}

A {\sf connection} on a left $A$-module $E$ is a $\C$-linear map $\nabla:E \to \Omega^1 A \otimes_A E$
such that
$$
\nabla(a.e) = \extd a \otimes e + a.\nabla e
$$
for all $a \in A$ and $e \in E$. Because $\Omega^1 A= \Omega^{1,0} A\oplus \Omega^{0,1}A$ we write $\nabla=\nabla^{1,0}+\nabla^{0,1}$ where $\nabla^{p,q}:E \to 
\Omega^{p,q}A \otimes_A E$.
We define higher covariant exterior derivatives 
$$
\nabla:\Omega^n A \otimes_A E \to \Omega^{n+1} A \otimes_A E \qquad \hbox{by}  \qquad
\nabla(\xi\tens e)\,:=\, \extd \xi \tens  e +(-1)^n\,\xi\wedge\nabla e
$$
for $\xi \in \Omega^nA$ and $e \in E$.
The {\sf curvature} of the connection $\nabla$ is 
$$
R:=\nabla \circ \nabla:E\to \Omega^2 A \otimes_A E.
$$
It can be shown that 
$$
\nabla \circ \nabla =\id \wedge R:\Omega^n A\otimes_A E\to  \Omega^{n+2}A \otimes_A E
$$
so, if the curvature vanishes, there is a complex $(\Omega^\hdot A \otimes_AE,\nabla)$ and an associated cohomology theory  \cite{bbsheaf}.

\medskip
Now we show how, in parallel with the Koszul-Malgrange Theorem,  
a connection induces a holomorphic structure if the $(0,2)$ part of its curvature vanishes. 
The vanishing of the holomorphic curvature is a weaker condition than the 
vanishing of the curvature for a standard covariant derivative, as shall be explained.

As before, $\pi^{p,q}:\Omega^{p+q}A \to \Omega^{p,q}A$ are the orthogonal projections.

\begin{propos}
Let $(E,\nabla)$ be a left $A$-module with a connection. Then the map
\begin{eqnarray*}
\nablaol\,:=\,(\pi^{0,1}\tens \id_E)\nabla :E\to\Omega^{0,1} A \otimes_A E
\end{eqnarray*}
is a $\pdol$-operator. 
Its holomorphic curvature in terms of the curvature $\nabla^2$ is given by the formula
\begin{equation}
\label{eq.hol.curv.formula}
\nablaol^2=(\pi^{0,2} \otimes \id_E) \nabla^2
%R^{\nablaol}_{\mathrm{hol}}=(\pi^{0,1}\wedge\pi^{0,1}\tens I_E)R
:E\to \Omega^{0,2}A \otimes_A E.
\end{equation}
Hence if the $(0,2)$ component of the curvature $\nabla^2$ vanishes, then $(E,\nablaol)$ is a holomorphic left $A$-module.
\end{propos}
\begin{pf}
First we check that $\nablaol$ is a $\pdol$-operator, i.e., that the left Leibniz rule holds.
Let $a \in A$ and $e \in E$. Because $\pi^{0,1}$ is a left $A$-module homomorphism
\begin{eqnarray*}
(\pi^{0,1}\tens \id_E)\nabla(a.e) &=& (\pi^{0,1}\tens \id_E)\,(\extd a\tens e+a.\nabla(e)) \cr
&=& \pi^{0,1}(\extd a)\tens e+a.(\pi^{0,1}\tens \id_E)\nabla(e)  \cr
&=& \bar\partial a\tens e + a.\nablaol(e)\ .
\end{eqnarray*}
Thus $\nablaol$ is a $\pdol$-operator.

If $\omega \in \Omega^{p+q}A$ we will write $\omega^{p,q}$ for $\pi^{p,q}(\omega)$, i.e., for the 
component of $\omega \in \Omega^{p,q}A$.

Let $e \in E$ and suppose that $\nabla e=\omega\tens f$ (or, more correctly, a sum of such terms).
Then $\nabla^2(e) = \nabla(\omega \otimes f)=  \extd \omega \otimes f - \omega \wedge \nabla f$ so
\begin{eqnarray*}
(\pi^{0,2} \otimes \id_E) \nabla^2(e) &=&\pi^{0,2}(\extd \omega)  \otimes f - (\pi^{0,2} \otimes \id_E)( \omega \wedge \nabla f ) \cr
 &=&\pi^{0,2}(\extd \omega^{1,0} + \extd \omega^{0,1})  \otimes f -  \omega^{0,1} \wedge  (\pi^{0,1} \otimes \id_E) ( \nabla f ) \cr
&=& \pdol \omega^{0,1}  \tens f -\omega^{0,1} \wedge\nablaol f \cr
&=& \nablaol(\omega^{0,1} \tens f) \cr
&=& \nablaol \big( (\pi^{0,1} \otimes \id_E)(\omega  \tens f)\big) \cr
&=& \nablaol \, \nablaol (e)
\end{eqnarray*}
so proving the formula in (\ref{eq.hol.curv.formula}).
The second line of the  calculation used the fact that if $\omega,\eta \in \Omega^1A$, then $\pi^{0,2}(\omega \wedge \eta)
= \omega^{0,1} \wedge \eta^{0,1}$; see (\ref{eq.wedge.mult}).

The last sentence of the proposition follows at once.
\end{pf}

\section{Dolbeault cohomology}

Let  $(\Omega^\hdot A,\extd,*)$ be a differential $*$-calculus with an integrable almost complex structure  $J$.

\subsection{}

\begin{defin} 
For each integer $p$ we call
\begin{equation}
\label{Db.cplx}
0 \longrightarrow \Omega^{p,0}A \stackrel{\bar\partial}\longrightarrow \Omega^{p,1}A 
 \stackrel{\bar\partial}\longrightarrow \Omega^{p,2}A \longrightarrow  \cdots
\end{equation}
the $p^{\th}$ {\sf Dolbeault complex}  and call its cohomology groups
$$
H^{p,q}_{\pdol}(A)
$$
the {\sf $(p,q)$-Dolbeault cohomology} of $A$.  
\end{defin}

The $q=0$ case of the classical theorem that the Dolbeault cohomology of a complex manifold $X$ computes the cohomology of the sheaf $\Omega^p_X$, i.e., 
$H^{p,q}(X) \cong H^q(X,\Omega^p_X)$, becomes a tautology here: $\Omega^p_{\sf hol}A$ is equal to $H_{\pdol}^{p,0}(A)$ because we {\it defined} $\Omega^p_{\sf hol}$ to be $\{\omega \in \Omega^{p,0} \; |
\; \pdol\omega =0\}$ in (\ref{defn.Omega.hol}).

 \subsection{}

More generally, if $(E,\nablaol)$ is a holomorphic $A$-module we define Dolbeault complexes 
\begin{equation}
\label{Db.cplx.E}
0 \longrightarrow \Omega^{p,0} A \otimes_AE \stackrel{\nablaol}\longrightarrow \Omega^{p,1}A \otimes_AE 
 \stackrel{\nablaol}\longrightarrow \Omega^{p,2}A \otimes_AE \longrightarrow  \cdots
\end{equation}
by 
$$
\nablaol(\xi \otimes e):= \pdol \xi \otimes e + (-1)^{|\xi|} \xi \wedge \nablaol e
$$
and Dolbeault cohomology groups  
$$
H^{p,q}(E,\nablaol).
$$

Let $M$ be a complex manifold and $\Omega_M$ its sheaf of holomorphic 1-forms, i.e., the 
holomorphic sections of the cotangent bundle. 
Then the sheaf of holomorphic $p$-forms on $M$ is $\Omega^p_M:=\bigwedge^p\Omega_M$.
If  $A$ is the $\C$-valued differentiable functions on $M$,  then 
$$
H^{p,q}_{\pdol}(A) \cong H^q(M,\Omega^p_M)
$$
 the $q^{\th}$ sheaf cohomology group  of the sheaf    $\Omega_M$; see e.g., \cite[Sect. 2.6]{Huy} and 
  \cite[Cor. 4.38 and Defn. 2.37]{V}. The dimensions $h^{p,q}$ of $H^{p,q}_{\pdol}(A)$ are the Hodge numbers. 
  
Similarly, $H^{p,q}(E,\nablaol) \cong H^q(M, E \otimes \Omega^p_M)$

\section{The Hodge to de Rham spectral sequence} \label{frolicher}

The  Hodge to de Rham spectral sequence for a complex manifold can be used to state 
some of the results of Hodge theory without assuming the underlying analysis and integrals
that underpin Hodge theory (see \cite{V}).

Let  $(\Omega^\hdot A,\extd,*)$ be a differential $*$-calculus with integrable almost complex structure  $J$.

The {\sf Hodge-Fr\"olicher}, or {\sf Hodge to de Rham},  spectral sequence is the spectral sequence 
associated to the double complex $\Omega^{\hdot,\hdot}A$ with its total differential $\extd=\pd+\pdol$.
The complex-valued  de Rham complex $(\Omega^\hdot,\extd)A$ has a decreasing filtration $F^p\Omega^\hdot A$ where
\begin{eqnarray*}
F^p\Omega^n A \,:=\,  \bigoplus_{p \le r} \Omega^{r,n-r}A\ .
\end{eqnarray*}
Since the quotient $F^p\Omega A/F^{p+1}\Omega A$ is the $p^{\rm th}$ Dolbeault complex
\begin{eqnarray}\label{ucfghcf}
\Omega^{p,0}A \stackrel{\bar\partial}\longrightarrow 
\Omega^{p,1}A \stackrel{\bar\partial}\longrightarrow
\Omega^{p,2}A \stackrel{\bar\partial}\longrightarrow\dots
\end{eqnarray}
the terms on the first page of the spectral sequence associated to this filtration are 
the Dolbeault cohomology groups, 
$$
E_1^{p,q}=H^{p,q}_{\pdol}(A),
$$
and the differential on $E_1^{\hdot,\hdot}$ is given by restricting 
$\partial:\Omega^{p,q}A \to \Omega^{p+1,q}A$
to the kernel of $\bar\partial$. This spectral sequence converges to
the cohomology of the total complex (i.e.\ $H^\hdot_{\mathrm{dR}}(A)$) in the sense that
\begin{eqnarray*}
E^{p,q}_\infty \,\cong\, \frac{F^pH^{p+q}_{\mathrm{dR}}(A)}{F^{p+1}H^{p+q}_{\mathrm{dR}}(A)}\ ,
\end{eqnarray*}
where $F^pH^{\hdot}_{\mathrm{dR}}(A)$ is the image of the map 
$H^\hdot(F^p\Omega A,\extd)\to H^\hdot_{\mathrm{dR}}(A)$ induced by the inclusion
$F^p\Omega A\to \Omega A$.

\begin{propos}
Suppose the integrable almost complex structure on $A$ satisfies the following condition:
\begin{quote}
\label{EJB.F.cond}
$\phantom{xxxxxx}$
for all $p,q\in\mathbb{N}$, the map $\wedge:\Omega^{0,q}A\otimes_A
\Omega^{p,0}A\to \Omega^{p,q}A$ 
\\
$\phantom{xxxxxx}$ 
 is an isomorphism   with inverse $\Theta^{p,q}$.
 \end{quote}
Then  each $\Omega^{p,0}A$ is a holomorphic left $A$-module with respect to the $\pdol$-operator
 \begin{eqnarray*}
\nablaol\,:=\, \Theta^{p,1}\bar\partial:\Omega^{p,0}A\to \Omega^{0,1} A \otimes_A \Omega^{p,0}A.
\end{eqnarray*}
Furthermore, the terms on the first page of the Fr\"olicher, or Hodge to de Rham,  spectral sequence are
$$
E^{p,q}_1 = H^{p,q}_{\pdol}(A) \cong H^q(\Omega^{p,0}A,\nablaol)
$$
and this spectral sequence converges to the de Rham cohomology of $A$.\footnote{Classically, if the conditions for Hodge theory are satisfied, this spectral sequence 
has all derivatives zero, and thus converges at the first page.} 
\end{propos}
\begin{pf} 
We need to check that the holomorphic curvature vanishes. By (\ref{eq.hol.curv}), the holomorphic curvature $\nablaol^2$ is the composition
$$
\UseComputerModernTips
  \xymatrix{ 
  \Omega^{p,0}A \ar[r]^>>>>>\nablaol  & \Omega^{0,1} A\otimes_A \Omega^{p,0}A \ar[rrr]^{\bar\partial\tens \id - \id \wedge \nablaol} &&&
\Omega^{0,2} A\otimes_A \Omega^{p,0}A\ .
}
$$
As $\wedge:\Omega^{0,2} A \otimes_A \Omega^{p,0}A \to \Omega^{p,2}A $ is an isomorphism, to show $\nablaol^2=0$ it suffices to show that the composition 
\begin{equation}
\label{eq.curv.3}
\UseComputerModernTips
  \xymatrix{ 
  \Omega^{p,0}A \ar[r]^>>>>>\nablaol  &  \Omega^{0,1}A \otimes_A \Omega^{p,0}A \ar[rrr]^{\bar\partial\tens \id - \id \wedge \nablaol} &&&
\Omega^{0,2} A \otimes_A \Omega^{p,0}A \ar[r]^>>>>>{\wedge} & \Omega^{p,2} A
}
\end{equation}
vanishes.

The associativity of the multiplication $\wedge$ implies that the diagram
$$
\UseComputerModernTips
  \xymatrix{ 
 \Omega^{0,1}A\otimes_A \Omega^{0,1}A \otimes_A \Omega^{p,0} A \ar[rr]^{\wedge \tens \id}
 \ar[d]_{\id \otimes \wedge} &&
\Omega^{0,2}A \otimes_A \Omega^{p,0} A \ar[d]^>>>>>{\wedge} 
\\
\Omega^{0,1}A \otimes_A \Omega^{p,1}  A\ar[rr]_{\wedge} && \Omega^{p,2} A
}
$$
commutes. By hypothesis, the vertical maps are isomorphisms with inverses $\id \otimes \Theta^{p,1}$ and $\Theta^{p,2}$ so the rectangle in the diagram
\begin{equation}
\label{big.diag}
\UseComputerModernTips
  \xymatrix{ 
  \Omega^{0,1}A \otimes_A \Omega^{p,0} A \ar[rr]^{\id \otimes \pdol} \ar[ddrr]_{\id \otimes \nablaol}
&& \Omega^{0,1}A \otimes_A \Omega^{p,1} A \ar[dd] | {\id \otimes \Theta^{p,1}} \ar[rr]^{\wedge} &&
   \Omega^{p,2} A\ar[dd]^{\Theta^{p,2}}
  \\
  \\
&& \Omega^{0,1} A \otimes_A \Omega^{0,1} A \otimes_A \Omega^{p,0} A \ar[rr]_{\wedge \tens \id}  &&
\Omega^{0,2}A \otimes_A \Omega^{p,0} A \ar[d]^{\wedge}
\\
&&&& \Omega^{p,2}A
}
\end{equation}
 commutes; the triangle commutes by the definition of $\nablaol$; by the commutativity of 
 (\ref{big.diag}),  
 $$
 \wedge \circ (\wedge \otimes \id) \circ (\id \otimes \nablaol) =
 \wedge \circ \Theta^{p,2} \circ \wedge \circ (\id \otimes \pdol) = \wedge \circ (\id \otimes \pdol).
 $$
 It follows that the composition in (\ref{eq.curv.3}) is equal to
 $$
 (\pdol \wedge \id - \id \wedge \pdol) \Theta^{p,1}\pdol.
 $$
 This is zero because $\bar\partial^2=0$. 
 
By the discussion prior to the lemma the Hodge to de Rham spectral sequence converges and $E_1^{p,q}$ is the Dolbeault cohomology group $H^{p,q}_{\pdol}(A)$. It therefore suffices to show that the 
$p^{\th}$ Dolbeault complex
\begin{equation}
  \label{diag1}
\Omega^{p,0} A \stackrel{\bar\partial}\longrightarrow 
\Omega^{p,1} A \stackrel{\bar\partial}\longrightarrow
\Omega^{p,2} A \stackrel{\bar\partial}\longrightarrow\dots
\end{equation}
is isomorphic to
\begin{equation}
 \label{diag2}
\Omega^{p,0} A \stackrel{\nablaol}\longrightarrow 
\Omega^{0,1} A \otimes \Omega^{p,0} A \stackrel{\nablaol}\longrightarrow
\Omega^{0,2} A \otimes \Omega^{p,0} A \stackrel{\nablaol}\longrightarrow\dots.
\end{equation}
  Let $\xi \otimes \eta \in \Omega^{0,q} A \otimes_A \Omega^{p,0}A$. The image of $\xi \otimes \eta$ in $\Omega^{p,q+1}A$ after going clockwise around the diagram
$$
\UseComputerModernTips
  \xymatrix{ 
 \Omega^{p,q} A \ar[rr]^{\pdol} && \Omega^{p,q+1}  A
  \\
  \Omega^{0,q} A \otimes_A \Omega^{p,0} A \ar[u]^{\wedge}   \ar[rr]_{\nablaol} &&
\Omega^{0,q+1}A \otimes_A \Omega^{p,0} A  \ar[u]_{\wedge}
}
$$
is $\pdol \xi \wedge \eta +(-1)^q \xi \wedge \pdol \eta$. Going counter-clockwise, $\xi \otimes \eta$ is first sent to 
$$\pdol \xi \otimes \eta +(-1)^q \xi \wedge \nablaol \eta = \pdol  \xi\otimes \eta +(-1)^q \xi \wedge \Theta^{p,1} \pdol \eta$$ which is then sent to  $\pdol \xi \wedge \eta +(-1)^q \xi \wedge \pdol \eta$. The rectangle therefore commutes and it follows that (\ref{diag1}) is isomorphic to  (\ref{diag2}).
\end{pf}

\begin{example}
The complex-valued  forms on a complex manifold are given locally 
in terms of the coordinates
 $z_i$ by 
wedges of $\extd z_i$ and $\extd\bar z_i$. To get the map $\Theta$, locally, just permute
all the $\extd\bar z_i$ to the left, introducing the appropriate power of -1,
and replace the $\wedge$ separating the $\extd\bar z_i$ from
the $\extd z_i$ by $\tens$. For example,
\begin{eqnarray*}
\extd z_1 \wedge \extd\bar z_2 &\longmapsto& -\,\extd\bar z_2 \tens \extd z_1\ ,\cr
\extd\bar z_1\wedge \extd z_1 \wedge \extd\bar z_2 &\longmapsto& -\,
\extd\bar z_1 \wedge \extd\bar z_2 \tens  \extd z_1\ ,\cr
\extd z_3\wedge \extd z_1 \wedge \extd\bar z_2 \wedge \extd\bar z_4 &\longmapsto&
\extd\bar z_2 \wedge \extd\bar z_4 \tens \extd z_3\wedge \extd z_1\ .
\end{eqnarray*}
\end{example}

 \section{Non-commutative complex differential geometry and non-commutative projective algebraic  geometry}
 \label{sect.nca+dg}

We begin this section with some speculations in \S\ref{sect.big.pic} about the existence of integrable almost complex structures on non-commutative complex projective varieties. 
A foundation for these speculations is provided by the examples in 
\S\ref{ssect.CPn-theta} and \S\ref{ssect.q.flag.vars}.

\subsection{Speculating about the big picture}
\label{sect.big.pic}
\label{ssect.defn.X}
Let $R$ be a $\C$-algebra having all the properties attributed to the ring $R$ 
at the beginning of \S\ref{ssect.ncag}. For brevity, we write $X=\Projnc(R)$. 
Following \cite{ATV1}, \cite{ATV2}, and \cite{AZ}, 
we write $\cO_X$ for the image of $R$ in $\Qcoh X$ and {\it define}
 $H^p(X,-) := \Ext^p_{\Qcoh X}(\cO_X,-)$.
 
This section speculates about non-commutative spaces $X$, $X_{\sf real}$,
${\bf X}$, ${\bf X}_{\sf real}$, and $Y$, and associated structures, that play roles analogous to the spaces in the right-hand rectangle below: 
$$
  \UseComputerModernTips
  \xymatrix{ 
  {\bf X}_{\sf real}-\{0\}  \ar[d]^{\R^+}\ar@{=}[r] & {\bf X}-\{0\} \ar[dd]^{\C^{\times}}
  &&
   \RR^{2n+2}-\{0\}  \ar[d]^{\R^{+}} \ar@{=}[r] & \C^{n+1}-\{0\} \ar[dd]^{\C^{\times}}
  \\
 \ar@/^3pc/[u]^{\sigma} Y \ar[d]^{U(1)}
 && &
 \ar@/^3pc/[u]^{\sigma} S^{2n+1} \ar[d]^{U(1)}
  \\
X_{\sf real} \ar@{=}[r]  & X
&&
  \C\PP^n \ar@{=}[r]  & \C\PP^n
  }
  $$
The left-hand columns represent the underlying real varieties of the complex 
varieties in the right-hand columns.  

The examples in sections \ref{ssect.CPn-theta} and \ref{ssect.q.flag.vars} suggest 
there should be finitely generated non-commutative
noetherian $\C$-algebras 
that correspond to the non-commutative spaces in the left-hand diagram just 
above and fit together as in the following picture:
$$
  \UseComputerModernTips
  
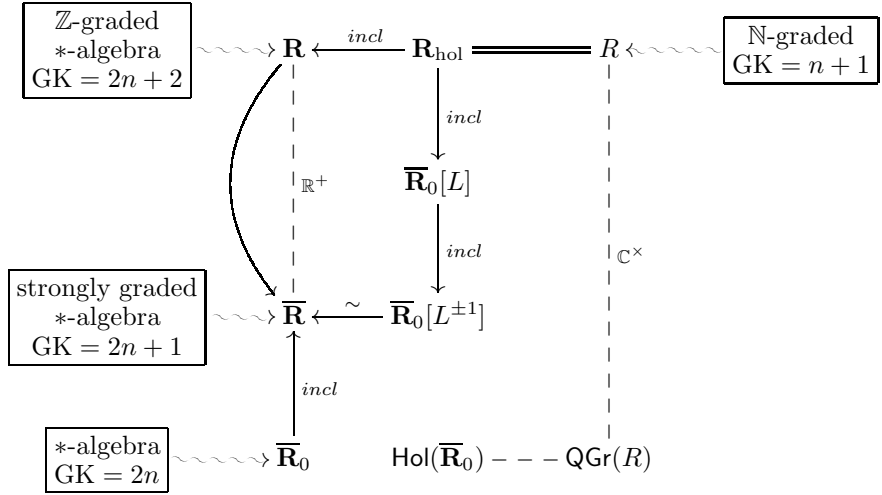
\xymatrix{ 
*+[F]{\txt{$\Z$-graded\\ $*$-algebra \\ $\GKdim =2n+2$}} \ar@{~>}[r] 
 & {\bf R}  \ar@/_2pc/[dd]   & \ar[l]_{incl} {\bf R}_{\rm hol} \ar[d]^{incl}& \ar@{=}[l]  R \ar@{--}[ddd]^{\C^{\times}} 
  &*+[F]{\txt{$\N$-graded\\ $\GKdim =n+1$}} \ar@{~>}[l]
     \\
     && \overline{\bf R}_0[L] \ar[d]^{incl}
     \\
*+[F]{\txt{strongly graded\\ $*$-algebra \\ $\GKdim =2n+1$}} \ar@{~>}[r]  & 
 \overline{\bf R}  \ar@{--}[uu]_{\RR^+} & \overline{\bf R}_0[L^{\pm 1}] \ar[l]_>>>>>>\sim
 && &
  \\
*+[F]{\txt{$*$-algebra \\ $\GKdim =2n$}} \ar@{~>}[r] 
& \overline{\bf R}_0 \ar[u]_{incl}    &{\sf Hol}(\overline{\bf R}_0) \ar@{--}[r]& \QGr(R)
  }
  $$
 The objects in this diagram  should have the following properties:   
  \begin{enumerate}
  \item{}
  ${\bf R}$ is a finitely generated noetherian $*$-algebra of Gelfand-Kirillov dimension $2n+2$
  and is a domain. It plays the role of $\C$-valued polynomial functions on a non-commutative
  real variety ${\bf X}_{\sf real}$ of real dimension $2n+2$.   
  \item 
 The $\C$-algebra ${\bf R}$ should have an integrable almost complex structure that 
 induces an integrable almost complex structure $\overline{\bf R}_0$. 
 \item{}
  ${\bf R}$ should have a $\Z$-grading that extends to $\Omega^\hdot{\bf R}$ and $({\bf R}_n)^*$
  should equal ${\bf R}_{-n}$.
      \item{}
  We will write $R$ for ${\bf R}_{\rm hol}$. 
  \item{}
  It is clear that $R$ is an $\N$-graded $\C$-algebra. We want $R_0=\C$ and $R$ should be generated by a finite set of elements $z_i \in R_1$. We also want $\GKdim(R)=n+1$ but it might be over optimistic to hope that  $R$ is noetherian. Let $\fm=R_{\ge 1}=(z_1,z_2,\ldots)$. 
  \item{}
  \label{genors}
  As a $\C$-algebra ${\bf R}$ should be generated by $R_1$ and its conjugate $R_1^*\subset {\bf R}_{-1}$.  
  Let ${\bf \fM}$ for the two-sided $*$-ideal of ${\bf R}$ generated by $R_1$. We want ${\bf R}/{\bf \fM} \cong \C$; if this is the case, then
    $\fm=R \cap {\bf \fM}$. Since $\C$ is the $*$-algebra of $\C$-valued functions on the point $\Spec(\R)$ the homomorphism ${\bf R} \to \C$ with kernel $\fM$ is the algebraic analogue of an inclusion 
  $\Spec(\R) \to {\bf X}_{\sf real}$; we think of the image of this inclusion as the vertex of the cone over $X$ and write 0 for that vertex. 
  \item
The ring $R$ plays two roles: that of the coordinate ring of the non-commutative complex affine variety ${\bf X}$, 
and that of the homogeneous coordinate ring of a non-commutative complex projective variety $X:=({\bf X}-\{0\})/\C^\times$. 
The objects in $\QGr(R)$ are ``quasi-coherent sheaves on $X$'' and objects in $\Mod(R)$  are ``quasi-coherent sheaves on 
${\bf X}$''; from this perspective, ${\bf X}$ is the cone over $X$.  The ``quasi-coherent sheaves'' on the punctured cone
${\bf X}-\{0\}$ are the objects in the quotient category $\Mod(R)/\Mod_0(R)$ where $\Mod_0(R)$ is the full subcategory of $\Mod(R)$
  consisting of the $R$-modules that are the sum of their finite dimensional submodules and those submodules are annihilated 
  by $\fm^n$ for $n \gg 0$. 
  \item{}
  Let $\Mod_0({\bf R})$ be the full subcategory of $\Mod({\bf R})$ consisting of the modules $M$ such that each element of $M$ is 
  annihilated by some power of $\fM$. Let us write $\Mod({\bf X}_{\sf real}-\{0\})$ for the quotient category $\Mod({\bf R})/\Mod_0({\bf R})$.
  \item 
  \label{norm.cond}
There should be a central element $c \in {\bf R}_0$ that is an $\R$-linear combination of the 
  elements  $z_iz_i^*$ and also an $\R$-linear combination of the 
  elements  $z_i^*z_i$.  The requirement that the coefficient of $z_iz_i^*$ belongs to $\R$ ensures that $c$, and hence $c-1$, is self-adjoint which 
  implies that ${\bf R}/(c-1)$ is a $*$-algebra. We think of $c$ as playing the the role of the ``norm-squared'' function.
  \item{}
  \label{norm=1}
  $\overline{\bf R}={\bf R}/(c-1)$ is the algebra of $\C$-valued polynomial functions on $Y$, 
  the norm=1 part of ${\bf X}_{\sf real}-\{0\}$, and ${\bf X}_{\sf real}-\{0\}$ behaves like $Y \times \RR^+$.
   \item{}
   Because $c-1$ is homogeneous, $\overline{\bf R}$ inherits a grading from ${\bf R}$. Its 
   degree-zero component   $\overline{\bf R}_0$ is the $U(1)$-invariant subalgebra,
       of $\overline{\bf R}$ and plays the role of $\C$-valued polynomial functions on $X_{\sf real}$,
  the non-commutative real variety underlying the non-commutative complex
  projective variety $X$ whose category of ``quasi-coherent sheaves'' is $\QGr(R)$.
  \item{}
  The solid arrows in the previous diagram represent $\C$-algebra homomorphisms and those labeled {\it incl}
  are injective.   The dashed lines $--$ should correspond to adjoint pairs of functors between certain abelian categories attached to the rings.
     \item{}
   ${\sf Hol}(\overline{\bf R}_0)$ is the category of holomorphic $\overline{\bf R}_0$-modules defined in section \ref{ssect.hol.mods} and should be related to $\QGr(R)$ in such 
a way that the cohomology groups $H^\hdot(E,\nablaol)$ for $(E,\nablaol) \in {\sf Hol}(\overline{\bf R}_0)$ coincide with appropriate cohomology groups $H^q(X,-)$ in $\QGr(R)$.
\item{}
\label{the.proj.L}
$L$ is an invertible $\overline{\bf R}_0$-bimodule, hence a rank one 
  projective $\overline{\bf R}_0$-module on both the left and the right,
  $$
  \overline{\bf R}_0[L^{\pm 1}] = \bigoplus_{m=-\infty}^\infty L^{\otimes m}
  \qquad \hbox{and} \qquad 
   \overline{\bf R}_0[L] = \bigoplus_{m=0}^\infty L^{\otimes m}.
   $$
\end{enumerate}

 \subsubsection{Complex line bundles on $X_{\sf real}$}
 
A closed subvariety $X$ of $\C\PP^n$ has a homogeneous coordinate ring $R$
 whose homogeneous components are (holomorphic) sections of complex line bundles on 
 $X_{\sf real}$. In the situation we are considering, invertible bimodules over $\overline{\bf R}_0$
 play the role of complex line bundles on $X_{\sf real}$. 
 
Each $\overline{\bf R}_m$ is an $\overline{\bf R}_0$-bimodule. 

From now on we will assume that $\overline{\bf R}$, and hence $\overline{\bf R}_0$, is a domain.
This is a mild hypothesis that roughly corresponds to $Y$ and $X_{\sf real}$ being irreducible varieties. 
  
Because $\overline{\bf R}_0$ has finite Gelfand-Kirillov dimension, it has a division ring of fractions, $D$ say.  The {\sf rank} of a finitely generated projective 
right $R$-module $P$ is $\dim_D (P \otimes_{\overline{\bf R}_0} D)$. 
In particular, if $P$ is isomorphic to a non-zero right ideal of $\overline{\bf R}_0$, then 
$\rank P=1$.  

If conditions (\ref{norm.cond}) and (\ref{norm=1}) above hold, then $1 \in \overline{\bf R}_1\overline{\bf R}_{-1}$ whence $\overline{\bf R}_0= \overline{\bf R}_1\overline{\bf R}_{-1}$. Likewise, 
$1 \in \overline{\bf R}_{-1}\overline{\bf R}_{1}$ and $\overline{\bf R}_0= \overline{\bf R}_{-1}\overline{\bf R}_{1}$. 

\begin{propos}
\label{prop.projective}
Suppose that $\overline{\bf R}$ is a domain and that conditions (\ref{genors}), (\ref{norm.cond}), 
and  (\ref{norm=1}) above hold. Then each $\overline{\bf R}_m$ is an invertible $\overline{\bf R}_0$-bimodule and a finitely generated rank one projective $\overline{\bf R}_0$-module on both the left and the right.
\end{propos}
\begin{pf}
For each integer $m$, the multiplication in $\overline{\bf R}$ gives an $\overline{\bf R}_0$-bimodule homomorphism
$$
\mu_m:\overline{\bf R}_m \otimes_{\overline{\bf R}_0} \overline{\bf R}_{-m}\to \overline{\bf R}_0.
$$
\underline{Claim:} $\mu_m$ is  surjective.  \underline{Proof}:
This is a triviality when $m=0$. 
The image of $\mu_m$ is a two-sided ideal of $R$ so it suffices to show 
that $1$ is in the image. The image of $\mu_1$ contains $z_jz_j^*$ for all $j$ and hence the
element $c$ in (\ref{norm.cond}); but $c=1$ in $\overline{\bf R}_0$ so $\mu_1$ is surjective. Similarly, $\mu_{-1}$ is surjective.

Suppose $m\ge 1$.
 Because $\overline{\bf R}$ is generated as a $\C$-algebra by 
 the image of $\{z_j^*,z_j \; | \; 0 \le j \le n\}$, 
$\overline{\bf R}_{\ge 0}$ is generated as an algebra by $\overline{\bf R}_0$
and $\overline{\bf R}_1$. 
Hence, if $m \ge 1$, $\overline{\bf R}_m=(\overline{\bf R}_1)^m$ by which we mean that $\overline{\bf R}_m$ is spanned  by the products of $m$ elements belonging to $\overline{\bf R}_1$. 
Likewise, $\overline{\bf R}_{-m}=(\overline{\bf R}_{-1})^m$. The image of $\mu_m$ is therefore 
$(\overline{\bf R}_1)^m(\overline{\bf R}_{-1})^m = \overline{\bf R}_1(\overline{\bf R}_{m-1} \overline{\bf R}_{1-m})\overline{\bf R}_1$. The truth of the claim for $m \ge 1$ now follows by induction. 

A similar argument shows $\mu_m$ is surjective when $m \le -1$. $\lozenge$ 

By the dual basis lemma, $\overline{\bf R}_m$ is a projective right $\overline{\bf R}_0$-module if the identity map is in the image of the map
$$
\Psi:\overline{\bf R}_m \otimes_{\overline{\bf R}_0} \Hom_{\overline{\bf R}_0}(\overline{\bf R}_m,(\overline{\bf R}_0)_{\overline{\bf R}_0}) \longrightarrow \End_{\overline{\bf R}_0}(\overline{\bf R}_m), \quad
\Psi(s \otimes \phi)(s'):=s\phi(s'),
$$
Given $a \in \overline{\bf R}_{-m}$ define $\phi_a \in \Hom_{\overline{\bf R}_0}(\overline{\bf R}_m,(\overline{\bf R}_0)_{\overline{\bf R}_0})$
by $\phi_a(s')=as'$. By the claim, there are elements  $s_i \in \overline{\bf R}_{m}$ and $a_i \in 
\overline{\bf R}_{-m}$ 
such that $\sum s_i a_i=1$. Hence 
$$
\Psi\Big(\sum_i s_i \otimes \phi_{a_i}\Big) = \id_{\overline{\bf R}_0}
$$
so $\overline{\bf R}_m$ is a projective right $\overline{\bf R}_0$-module.  A similar argument shows that $\overline{\bf R}_m$ is a projective left $\overline{\bf R}_0$-module. 

%We can interpret $\mu_m$ as giving an injective 
%homomorphism $L_m  \to L_{-m}^\vee := \Hom_{R}(L_m,R)$.
 
Because $\overline{\bf R}$ is a domain, left multiplication by a fixed non-zero element in $\overline{\bf R}_{-m}$  is an injective homomorphism $\overline{\bf R}_m \to \overline{\bf R}_0$ of right $\overline{\bf R}_0$-modules. Hence $\rank \overline{\bf R}_m=1$ on the right. As similar argument shows that $\overline{\bf R}_m$ has rank 1 as a left $\overline{\bf R}_0$-module.

The invertibility of $\overline{\bf R}_m$ as an  $\overline{\bf R}_0$-bimodule is part (3) of the next corollary.
\end{pf}

\begin{cor}
Suppose that $\overline{\bf R}$ is a domain and that conditions (\ref{genors}), (\ref{norm.cond}), 
and  (\ref{norm=1}) above hold. Then
\begin{enumerate}
  \item 
  $\overline{\bf R}$ is a strongly graded $\C$-algebra. 
  \item 
  There is an equivalence of categories ${\Mod}(\overline{\bf R}_0)  \equiv   {\Gr}(\overline{\bf R})$.
  \item 
  The multiplication map 
$
\overline{\bf R}_k \otimes_{\overline{\bf R}_0} \overline{\bf R}_m \to \overline{\bf R}_{k+m}
$
is an isomorphism of $\overline{\bf R}_0$-bimodules for all $k,m \in \Z$.
\item
Condition (\ref{the.proj.L}) holds.
\end{enumerate}
 \end{cor}
\begin{pf}
(1)
By definition, a $\Z$-graded ring $A$ is strongly graded if $A_iA_{-i}=A_0$ for all $i$. The proof of Proposition \ref{prop.projective}
shows that $1 \in \overline{\bf R}_{-m}\overline{\bf R}_{m}$ for all $m \in \Z$ so $\overline{\bf R}$ is strongly graded.

(2), (3), and (4), are standard facts about strongly graded rings.
\end{pf}

Given Proposition \ref{prop.projective}, the philosophy of non-commutative geometry says 
that elements of $\overline{\bf R}_m$ are smooth sections of a complex line bundle  on 
the real algebraic variety $X_{\sf real}$. 
After we give $X_{\sf real}$, equivalently $\overline{\bf R}_0$, an integrable
almost complex structure there should be a $\pdol$-operator $\nablaol$ 
on each $\overline{\bf R}_m$ such that $R_m$ consists of the elements of  $\overline{\bf R}_m$
on which $\nablaol$ vanishes.

\subsubsection{The (punctured) affine cone over $X$}
As in algebraic geometry, there is a non-commutative complex affine variety ${\bf X}$ that is a cone over $X$ defined implicitly by saying that the category of quasi-coherent sheaves on ${\bf X}$
is $\Mod R$. We then think of $R$ as a ring of  ``holomorphic functions'' on ${\bf X}$. 

The Artin-Schelter regularity conditions require $R$ to have good homological properties that
one should probably think of as saying ${\bf X}$ is smooth, so not quite what most people would
think of when they say ``cone''. If necessary, one could weaken the hypothesis on $R$
 by allowing it to be a quotient of an Artin-Schelter regular ring by a central regular sequence
such that the vertex of the cone is the only singular point, i.e., the ``deviation'' of $R$ from being Artin-Schelter regular is ``concentrated'' at the module $R/R_{\ge 1}$. Formally, one would ask that
the quotient category $(\Mod R)/{\sf S}$, where ${\sf S}$ is the smallest localizing subcategory containing 
$R/R_{\ge 1}$, has finite global dimension, i.e., $\Ext^{2n+1}$ vanishes on $(\Mod R)/{\sf S}$.
In such a  case it is more appropriate to  think of $R$ as a ring of ``holomorphic functions''
on ${\bf X}-\{0\}$.

\subsubsection{The action of $\C^\times$ on ${\bf X}$}
The homogeneous components of $R$ are the eigenspaces for an action of $\C^\times$ as automorphisms of $R$. If we think of this as being inherited from an action of $\C^\times$ on ${\bf X}$, then $\Gr R$ is, in effect, the category of $\C^\times$-equivariant 
(i.e., constant along orbits) 
quasi-coherent modules on ${\bf X}$. The {\sf vertex} $0 \in {\bf X}$ is, by definition, the unique point fixed by 
$\C^\times$; formally, $\{0\}=\Spec(R/R_{\ge 1})$. By quotienting out $\Fdim R$ we are removing those modules supported at the vertex. Thus $X$ behaves as the orbit space $({\bf X}-\{0\})/\C^\times$. 

\subsubsection{Viewing ${\bf X}$ as a non-commutative real algebraic variety with an almost complex structure}
\label{sect.real.R}

Viewed through the lens of the ring $R$, ${\bf X}$ is a non-commutative 
{\it complex} variety of complex dimension 
$n+1$. If we wish to treat ${\bf X}$ as a {\it real} algebraic variety of real dimension $2(n+1)$ 
we need a $*$-algebra ${\bf R}$ of GK-dimension $2n+2$ that will play the role of $\C$-valued differentiable functions on ${\bf X}_{\sf real}-\{0\}$.

An appropriate almost complex structure on ${\bf X}_{\sf real}$  would be an almost 
complex structure $(\Omega^\hdot {\bf R},\extd,*,J)$ such that 
\begin{enumerate}
  \item 
  ${\bf R}_{\sf hol}=R$,
  \item 
  there is a $\Z$-grading on ${\bf R}$ that is compatible with that on $R$, and
  \item
  $({\bf R}_n)^*= {\bf R}_{-n}$.
\end{enumerate}

Because we want to retain the underlying real structure on $X$ we use the factorization 
$\C^\times =\R^+ \times U(1)$ to pass down the left-hand sides of the diagrams at the beginning of 
section \ref{sect.big.pic}. \S\S\ref{ssect.R+bun} and \ref{ssect.U1.action} address this.

\subsubsection{${\bf X}_{\sf real}-\{0\}$ as a trivial $\RR^+$-bundle over $Y$} 
\label{ssect.R+bun}

Let $S^{2n+1}$ be the sphere of radius 1 in $\RR^{2n+2}$ centered at the origin. 
Then $\RR^{2n+2}-\{0\}$ is diffeomorphic to $\RR^+ \times S^{2n+1}$. 
A non-commutative analogue of the inclusion $S^{2n+1} \to \RR^{2n+2}$ is a surjective homomorphism 
$$
{\bf R} \to \overline{\bf R}:=\frac{\bf R}{(c-1)}
$$ 
where $c-1$ is a self-adjoint central element in ${\bf R}_0$. 
Since $\overline{\bf R}$ is to play the role of $\C$-valued polynomials on a non-commutative real algebraic 
variety $Y$ of dimension $2n+1$ we want the GK-dimension of $\overline{\bf R}$ to be $2n+1$. By \cite[Prop. 3.15]{KL}
$\GKdim(\overline{\bf R}) \le 2n+1$. 

The homomorphism ${\bf R} \to \overline{\bf R}$  induces a closed immersion $\s:Y \to {\bf X}_{\sf real}-\{0\}$ and there should be 
a non-commutative map $\pi:{\bf X}_{\sf real}-\{0\} \to Y$ such that $\pi\s=\id_Y$. 
The existence of such $\pi$ and $\s$ can be phrased in terms of the associated 
module categories by using the ideas and language in \cite{VdB} and \cite{Sm}. 
A {\sf closed immersion} is a triple of functors $(\s^*,\s_*,\s^!)$ such that 
$\s_*:\Mod\big(\!{\bf R}/(c-1)\big) \to \Mod({\bf X}_{\sf real}-\{0\})$ is fully faithful, the essential image of $\s_*$ is closed under 
subquotients, and $\s_*$ has a left adjoint $\s^*$ and a right adjoint $\s_*$. 

For example, the functor $\b_*:\Mod\big({\bf R}/(c-1)\big) \to \Mod({\bf R})$ that sends an ${\bf R}/(c-1)$ module $M$ to $M$ viewed as 
an ${\bf R}$-module determines a closed immersion $\b:Y \to {\bf X}_{\sf real}$. Similarly, The homomorphism ${\bf R} \to {\bf R}/\fM \cong \C$ corresponds to a closed immersion $\a:\{0\} \to {\bf X}_{\sf real}$. 

The localization functor $j^*:\Mod({\bf R}) \to \Mod({\bf X}_{\sf real}-\{0\})$ and its right adjoint $j_*$ implicitly define an open immersion $j:{\bf X}_{\sf real} -\{0\} \to {\bf X}$. 

\begin{lemma}
Let $\b_*$ be the forgetful functor associated to the homomorphism ${\bf R} \to {\bf R}/(c-1)$, and let $\b^*$ and 
$\b^!$ the usual left and right adjoints to $\b_*$. 
The functors $\s^*:=\b^*j_*$, $\s_*:=j^*\b_*$, and $\s^!:=\b^!j_*$, define a closed immersion 
$\s:Y \to {\bf X}_{\sf real}-\{0\}$ such that  the diagram
$$
\UseComputerModernTips
\xymatrix{
{\bf X}_{\sf real}-\{0\}  \ar[r]^j & {\bf X}_{\sf real} & \ar[l]_\a \{0\}
\\
Y   \ar@/^3pc/[u]^{\sigma} \ar[ur]_{\b} 
}
$$
commutes. Furthermore, $\s^*\s_* \cong \id$ and $\s^!\s_* \cong \id$.
\end{lemma}
\begin{pf}
Because $\b_*$ is left adjoint to $\b^!$ and $j^*$ is left adjoint to $j_*$, $j^*\b_*$ is left adjoint to $\b^!j_*$.

To key step in showing $\b^*j_*$ is left adjoint to $\s_*:=j^*\b_*$ is
to show that $\Hom_{\bf R}(M,N)=\Ext^1_{\bf R}(M,N)=0$ for all ${\bf R}/(c-1)$-modules 
$N$ and all ${\bf R}/\fM$-modules $M$  (see \cite[Prop. 7.1]{Sm} for example). Since every ${\bf R}/\fM$-module is a direct sum of copies of ${\bf R}/\fM$ it suffices to show that $\Hom_{\bf R}({\bf R}/\fM,N)=\Ext^1_{\bf R}({\bf R}/\fM,N)=0$ for all ${\bf R}/(c-1)$-modules $N$. 

Since $c$ belongs to $\fM$, $R/\fM$ is not annihilated by $c-1$, whence $\Hom_{\bf R}({\bf R}/\fM,N)=0$. 

Now consider an exact sequence $0 \to N \to E \to {\bf R}/\fM \to 0$ of left
${\bf R}$-modules such that $(c-1)N=0$. Since $c$ is central multiplication by $c-1$ gives an ${\bf R}$-module homomorphism $E \to (c-1)E$. Hence $(c-1)E \cong {\bf R}/\fM$. Since $(c-1)N=0$, 
$E=N \oplus (c-1)E$ and this leads to a splitting of the sequence $0 \to N \to E \to {\bf R}/\fM \to 0$. 

We now show $\s^*$ is left adjoint to $\s_*$.  Let $N$ be an ${\bf R}$-module. It is a standard fact \cite{Gabriel} 
about localizing subcatgeories that there is an exact sequence 
$$
0 \to \varinjlim_n \Hom_{\bf R}({\bf R}/\fM^n,N) \to N \to j_*j^*N \to  \varinjlim_n \Ext^1_{\bf R}({\bf R}/\fM^n,N) \to 0.
$$
From this sequence and what we proved above we see that the natural map
$N \to j_*j^*N$ is an isomorphism for all ${\bf R}/(c-1)$-modules $N$. Hence $j_*j^*\b_* \cong \b_*$. It is also a standard fact about localizations
that $j^*j_* \cong \id$.  Let $N$ be an ${\bf R}/(c-1)$-module. Then
$$
 \Hom(\b^*j_*M,N) \cong    \Hom(j_*M,j_*j^*\b_*N) \cong  \Hom(j^*j_*M,j^*\b_*N)
\cong \Hom(M,j^*\b_*N).
$$
Thus $\s^*$ is left adjoint to $\s_*$. It is easy to check that $\s^*\s_* \cong \id$ and $\s^!\s_* \cong \id$. 
\end{pf}

For this to be compatible with the principal $\C^\times$-bundle on the 
right-hand sides of the diagrams at the beginning of 
section \ref{sect.big.pic} it should be the case that ${\bf R}_m=L^{\otimes m}$. 

As mentioned before the previous result we would like a non-commutative map 
$\pi:{\bf X}_{\sf real}-\{0\} \to Y$ such that $\pi\s=\id_Y$. Since $\s^*\s_* \cong \id$ one might define $\pi_*= \s^*$. However, we do not know whether this is appropriate since it doesn't appear that such a $\pi_*$ has a left adjoint.

\subsubsection{The action of $\U(1)$ on ${Y}$}
\label{ssect.U1.action}
There is an action of $U(1)$ as $*$-algebra automorphisms of ${\bf R}$ and ${\bf R}/(c-1)$ given by $\xi \cdot r = \xi^m r$ if
$r \in L^{\otimes m}$. 
If $x \in {\bf R}_m$ and $\xi \in U(1)$, then $\xi\cdot x^*=(\xi \cdot x)^* = (\xi^m x)^*=\xiol^m x^* = \xi^{-m} x^*$ so $x^* \in {\bf R}_{-m}$.   

\subsubsection{Analogues of the sheaves $\Omega_X^p$}
\label{defn.OmegaX}

Even for those non-commutative projective varieties $X$ that are well understood
and have excellent homological properties one knows no objects in 
$\Qcoh(X)$ that deserve being denoted  
$\Omega^p_X$. 
There is no effective general theory of exterior powers in non-commutative ring theory so
it is unwise to seek objects   $\Omega^p_X$ that are exterior powers of $\Omega^1_X$. 

An alternative is to ask that $R$, which is connected graded, 
be a twisted Calabi-Yau algebra of dimension
$n$, i.e., the projective dimension of $R$ as an $R$-bimodule is  $n$  and the shifted 
Hochschild cohomology $HH^{\hdot}(R, R \otimes R)[n]$ is an invertible $R$-bimodule. 
Over a connected graded algebra the only invertible graded bimodules are of the form 
${}^1R^\sigma$ for some graded algebra automorphism $\sigma$ of $R$.
Because $R$ is connected graded  it has a unique minimal projective resolution as an
$R$-bimodule. The wished for  $\Omega^p_X$  might be defined in terms of this resolution. 
Beilinson's resolution of the diagonal is the motivation for this approach. 
 
A test as to whether these ideas  are promising is to ask if, after proceeding as above,
for some $X$'s that we understand, $H^{p,q}_{\pdol}({\bf R})$ is isomorphic to 
$H^q(X,\Omega^p_X)$ where the latter is defined as $\Ext^q_{\Qcoh X}(\cO_X,\Omega^p_X)$.

\subsubsection{Testing these ideas on non-commutative analogues of $\PP^n_{\C}$}

Let $R$ be an Artin-Schelter regular algebra \cite{AS} such that the minimal resolution
of the trivial module $\C=R/R_{\ge 1}$ is 
$$
0 \to R(-n) \to R(-n+1)^{n} \to \cdots \to R(-1)^{n} \to R \to \C \to 0.
$$
For such an $R$ it is generally accepted that $X$ is a good non-commutative analogue of 
$\PP^n_\C$. 
The image in $\Qcoh(X)$ of the resolution is an exact sequence 
$$
0 \to \CO(-n) \to \CO(-n+1)^{n} \to \cdots \to \CO(-1)^{n} \to  \CO_X \to 0.
$$
When $X$ {\it is} $\PP^n_\C$,  
\begin{equation}
\label{eq.sheaf.Omega^p.inQGr.R}
\ker \Big( \CO(-p)^{\binom{n}{p}} \to  \CO(-p+1)^{\binom{n}{p-1}}\Big) = \Omega^p_{\rm hol} 
\end{equation}
and Bott's Vanishing Theorem (see, e.g., \cite[Thm. 7.2.3]{CMP}) 
says $H^q(\PP^{n-1}_\C,\Omega^p_{\rm hol}(k))$ vanishes except for
\begin{enumerate}
  \item 
  $p=q$ and $k=0$;
  \item 
  $q=0$ and $k>p$;
  \item 
  $q=n-1$ and $k<-n+1+p$.
\end{enumerate}

The object $\Omega^p_{\rm hol}$ in  (\ref{eq.sheaf.Omega^p.inQGr.R})
is not the same as the object called $\Omega^p_{\sf hol}$ in  (\ref{defn.Omega.hol});
the latter was defined to be
 $\{\omega \in \Omega^{p,0}{\bf R} \; | \; \pdol\omega =0\}$. However, if $\{\omega \in \Omega^{p,0}{\bf R} \; | \; \pdol\omega =0\}$ is a graded $R$-module it would have an image in the quotient
category $\Qcoh(X)$ and we could ask if that image is isomorphic to the kernel in (\ref{eq.sheaf.Omega^p.inQGr.R}). If so, Bott's Vanishing Theorem would hold for $X$: 
because $R$ is Gorenstein, which is part of the definition of Artin-Schelter regularity, 
the cohomology  for $\cO_X(k)$ is the same as that for $\cO_{\PP^{n-1}}(k)$. 

For $\Omega^p_{\sf hol}$, as defined in  (\ref{defn.Omega.hol}),
to be a graded $R$-module  the $\C^\times$-action on 
 ${\bf R}$ must extend to a  suitable $\C^\times$-action on  $\Omega^{\hdot}{\bf R}$.

\subsection{An example: $\theta$-deformed $\mathbb{C}^{n+1}$ and $\mathbb{C}\PP^{n}$}
\label{ssect.CPn-theta}

There is a large physics literature on the subject of $\theta$-deformed spaces, or Moyal products. This has been used in quantisation, gauge theory, string theory, cosmology and integrable systems. 
We will quote the relations for the algebra and differential calculus from \cite{MarPenPanOys}.

\subsubsection{The noncommutative $(2n+2)$-plane $\mathbb{R}^{2n+2}_\theta$}

Let $\theta$ be a real skew-symmetric $(n+1)\times (n+1)$ matrix, and set $\lambda^{\mu\nu}=e^{\mathrm{i}\theta_{\mu\nu}}$. 

The algebra of complex-valued polynomial functions on $\mathbb{R}^{2n+2}_\theta$ is the 
$*$-algebra $\C\big[\mathbb{R}^{2n+2}_\theta\big]$
%$\C_{\mathrm{alg}}(\mathbb{R}^{2n}_\theta)$
 generated by 
$\{z^\mu, \bar{z}^\mu \; | \; 0\le\mu\le n\}$ with relations
\begin{eqnarray}
\label{eq.z.relns}
z^\mu z^\nu=\lambda^{\mu\nu}z^\nu z^\mu\ ,\quad \bar{z}^\mu \bar{z}^\nu=\lambda^{\mu\nu}\bar{z}^\nu \bar{z}^\mu\ ,\quad
\bar{z}^\mu z^\nu = \lambda^{\nu \mu} z^\nu \bar{z}^\mu\,
\end{eqnarray}
and $(z^\mu)^*=\bar{z}^\mu$. 
%$C_{\mathrm{alg}}(\mathbb{R}^{2n}_\theta)$ is
The elements 
$$
(z^0)^{i_0}\dots (z^n)^{i_n}\,(\bar{z}^0)^{j_0}\dots (\bar{z}^n)^{j_n}\ ,\quad (i_0,\dots,i_n,j_0,\dots,j_n)\in\mathbb{N}^{2n+2}\ 
$$
are a $\mathbb{C}$-vector space basis for $\C\big[\mathbb{R}^{2n+2}_\theta\big]$.
 
There is an action of $\C^\times$ as algebra automorphisms of 
$\C\big[\mathbb{R}^{2n+2}_\theta\big]$
defined by 
%$C_{\mathrm{alg}}(\mathbb{R}^{2n}_\theta)$ given by 
$\xi\cdot z^\mu=\xi\, z^\mu$ and
$\xi\cdot \bar{z}^\mu=\xiol\, \bar{z}^\mu$. Its subgroup $U(1)$ acts as 
$*$-algebra automorphisms of $\C\big[\mathbb{R}^{2n+2}_\theta\big]$.
The $\C^\times$-action induces a $\Z$-grading on $\C\big[\mathbb{R}^{2n+2}_\theta\big]$ with 
$\deg z^\mu=1$ and $\deg \bar{z}^\mu=-1$. 

Since $\C\big[\mathbb{R}^{2n+2}_\theta\big]$
%$C_{\mathrm{alg}}(\mathbb{R}^{2n}_\theta)$ 
is a quadratic algebra having a PBW basis it is a Koszul algebra by  Priddy's Theorem \cite[Thm. 5.3]{Priddy}.

\subsubsection{A differential $*$-calculus on $\mathbb{R}^{2n+2}_\theta$}

A $*$-differential calculus $(\Omega^\hdot,\extd,*)$ on $\C\big[\mathbb{R}^{2n+2}_\theta\big]$ is defined in \cite{MarPenPanOys}: $\Omega^\hdot \C\big[\mathbb{R}^{2n+2}_\theta\big]$ has generators $z^\mu,\bar{z}^\mu \in \Omega^0$ and $\extd z^\mu, \extd\bar{z}^\mu \in \Omega^1$ with relations  (\ref{eq.z.relns}) and
\begin{eqnarray}
z^\mu \extd z^\nu & = & \lambda^{\mu\nu}\extd z^\nu z^\mu\ , \cr
\bar{z}^\mu \extd\bar{z}^\nu & = & \lambda^{\mu\nu} \extd\bar{z}^\nu \bar{z}^\mu\ ,\cr
\bar{z}^\mu \extd z^\nu & = & \lambda^{\nu\mu}\extd z^\nu \bar{z}^\mu\ , \cr
z^\mu \extd\bar{z}^\nu & = & \lambda^{\nu\mu} \extd\bar{z}^\nu z^\mu\ ,
\end{eqnarray}
and  
\begin{eqnarray}
\extd z^\mu \wedge \extd z^\nu+\lambda^{\mu\nu}\extd z^\nu \wedge  \extd z^\mu  &= & 0\ , \cr
\extd\bar{z}^\mu \wedge  \extd\bar{z}^\nu+\lambda^{\mu\nu} \extd\bar{z}^\nu \wedge  \extd\bar{z}^\mu  &= & 0\ ,\cr
\extd\bar{z}^\mu \wedge  \extd z^\nu+\lambda^{\nu\mu}\extd z^\nu \wedge  \extd\bar{z}^\mu  & = & 0,
\end{eqnarray}
and $*$-structure defined by
$$
(z^\mu)^*=\bar{z}^\mu
\qquad \hbox{and} \qquad (\extd z^\mu)^*=\extd \bar{z}^\mu
$$
and extended to $\Omega^\hdot\C\big[\mathbb{R}^{2n+2}_\theta\big]$ by $(\omega\eta)^*:=(-1)^{mn}\eta^*\omega^*$ for $\omega \in \Omega^m\C\big[\mathbb{R}^{2n+2}_\theta\big]$ and $\eta \in \Omega^n\C\big[\mathbb{R}^{2n+2}_\theta\big]$.

The $\C^\times$-action extends to $\Omega^\hdot\C\big[\mathbb{R}^{2n+2}_\theta\big]$ by $\xi\cdot \extd z^\mu=\xi\, \extd z^\mu$ and
$\xi\cdot \extd\bar{z}^\mu=\xiol \, \extd\bar{z}^\mu$.

\subsubsection{An integrable almost complex structure on $\mathbb{R}^{2n+2}_\theta$} \label{kjhsdgvjkh}

There is a unique degree-preserving map $J:\Omega^\hdot\C\big[\mathbb{R}^{2n+2}_\theta\big] \to \Omega^\hdot\C\big[\mathbb{R}^{2n+2}_\theta\big]$ that is
\begin{enumerate}
  \item 
an $A$-bimodule endomorphism, and
  \item 
a derivation, and satisfies  
  \item 
  $J(\extd z^\mu):=\mathrm{i}\,\extd z^\mu$ and $J(\extd\bar{z}^\mu):=-\mathrm{i}\,\extd\bar{z}^\mu$.
\end{enumerate} 
In other words, $J$ is an almost complex structure on $\mathbb{R}^{2n+2}_\theta$.  We 
therefore have spaces $\Omega^{p,q}\C\big[\mathbb{R}^{2n+2}_\theta\big]$ and derivations $\pd$ and $\pdol$ defined as in \S\ref{ssect.Omega.pq} and \S\ref{ssect.d.dbar}.

The $\C^\times$-action extends with
$\xi\cdot\omega=\xi^{p}\, \bar{\xi}^{q}\,\omega$ for $\omega\in\Omega^{p,q}\C\big[\mathbb{R}^{2n+2}_\theta\big]$. 

\begin{lemma}
$J$ is integrable.
\end{lemma}
\begin{pf}
Every element in $\Omega^{1,0}\C\big[\mathbb{R}^{2n+2}_\theta\big]$ can be written in the form
$f_\mu.\extd z^\mu$ for some $f_\mu\in \C\big[\mathbb{R}^{2n+2}_\theta\big]$. Applying $\extd$ to this gives
$\extd f_\mu\wedge\extd z^\mu$, which is in $\Omega^{2,0}\C\big[\mathbb{R}^{2n+2}_\theta\big] \oplus \Omega^{1,1}\C\big[\mathbb{R}^{2n+2}_\theta\big]$. 
Hence $\extd$ satisfies condition (4) in Lemma~\ref{lem.int}.
% namely $\extd \Omega^{1,0} \subset \Omega^{2,0} \oplus \Omega^{1,1}$.
\end{pf}

The subalgebra $\C\big[\mathbb{R}^{2n+2}_\theta\big]_{\sf hol}$ is  generated by the $z^\mu$, $0 \le \mu \le n$. We denote it by $\cO(\C^{n+1}_\theta)$ and think of it as the coordinate ring of a 
non-commutative affine variety $\C^{n+1}_\theta$. The algebra $\cO(\C^{n+1}_\theta)$ has all the properties attributed to the ring $R$ at the beginning of \S\ref{ssect.ncag}.
It is well-known that $\Projnc \cO(\C^{n+1}_\theta)$ behaves very much like $\C\PP^n$.

\subsubsection{The noncommutative sphere $S^{2n+1}_\theta$}

For each $\mu$, the element $\bar{z}^\mu z^\mu=z^\mu\bar{z}^\mu$ is central.
We define the $\C$-valued polynomial functions on $S^{2n+1}_\theta$ to be
$$
\C\big[S^{2n+1}_\theta\big]:= \frac{\C\big[\mathbb{R}^{2n+2}_\theta\big]}{\big(\sum_\mu z^\mu\bar{z}^\mu -1\big)}.
$$
Since $\sum_\mu z^\mu\bar{z}^\mu -1$ is self-adjoint $\C\big[S^{2n+1}_\theta\big]$ inherits
a $*$-algebra structure from $\C\big[\mathbb{R}^{2n+2}_\theta\big]$.
Since $\sum_\mu z^\mu\bar{z}^\mu -1$ is invariant under the $\C^\times$-action 
the unit circle $U(1)$ acts as $*$-algebra automorphisms of $\C\big[S^{2n+1}_\theta\big]$. 

Differentiating the relation $\sum_\mu z^\mu\bar{z}^\mu=1$ gives 
$\sum_\mu \extd z^\mu\bar{z}^\mu+\sum_\mu z^\mu\extd\bar{z}^\mu=0$, so we define 
$\Omega^\hdot \C\big[S^{2n+1}_\theta\big]$ to be the quotient of  $\Omega^\hdot \C\big[\R^{2n+2}_\theta\big]$ by the relations $\sum_\mu z^\mu\bar{z}^\mu=1$ and 
$\sum_\mu \extd z^\mu\bar{z}^\mu+\sum_\mu z^\mu\extd\bar{z}^\mu=0$.

\subsubsection{$\C\PP^{n}_\theta$}
We define  $\C\big[\mathbb{CP}^{n}_\theta\big]$ to be the $U(1)$-invariant subalgebra 
of $\C\big[S^{2n+1}_\theta\big]$. Since $U(1)$ acts as $*$-algebra automorphisms, 
$\C\big[\mathbb{CP}^{n}_\theta\big]$ is a $*$-algebra. 

We define $\Omega^\hdot \C\big[\mathbb{CP}^{n}_\theta\big]$ to be the $U(1)$-invariant subalgebra of $\Omega^\hdot \C\big[S^{2n+1}_\theta\big]$, with the additional relations
$\sum_\mu \extd z^\mu\bar{z}^\mu=0$ and $\sum_\mu  z^\mu\extd\bar{z}^\mu=0$.

These relations $\sum_\mu \extd z^\mu\bar{z}^\mu=0$ in $\Omega^{1,0} \C\big[\mathbb{CP}^{n}_\theta\big]$ and $\sum_\mu  z^\mu\extd\bar{z}^\mu=0$ in  $\Omega^{0,1} \C\big[\mathbb{CP}^{n}_\theta\big]$ correspond to the following projection matrices for finitely generated projective modules:
$$
P_{\mu\nu}\,=\, \delta_{\mu\nu}-z^\mu\bar{z}^\nu\ ,\quad Q_{\mu\nu}\,=\, \delta_{\mu\nu}-\bar{z}^\mu z^\nu\ .
$$
The complex structure is given by $J$ from Section~\ref{kjhsdgvjkh}, i.e.\ 
$J(\extd z^\mu):=\mathrm{i}\,\extd z^\mu$ and $J(\extd\bar{z}^\mu):=-\mathrm{i}\,\extd\bar{z}^\mu$. A brief check shows that this is consistent with the additional relations given here. 

\subsubsection{Complex line bundles on $\C\PP^{n}_\theta$}

 The $U(1)$-action on $\C\big[S^{2n+1}_\theta\big]$ corresponds to a $\Z$-grading on $\C\big[S^{2n+1}_\theta\big]$  and  $\C\big[\mathbb{CP}^{n}_\theta\big]$ is its degree-zero component.
The homogeneous components of $\C\big[S^{2n+1}_\theta\big]$  are therefore
$\C\big[\mathbb{CP}^{n}_\theta\big]$-bimodules. 

They are, in fact, invertible $\C\big[\mathbb{CP}^{n}_\theta\big]$-bimodules, hence 
finitely generated rank one projective $\C\big[\mathbb{CP}^{n}_\theta\big]$-modules on both the right and left so should be thought of as polynomial sections of complex line bundles on 
$\mathbb{CP}^{n}_\theta$. The holomorphic sections should be thought of as the global sections of 
invertible $\cO_{\mathbb{CP}^{n}_\theta}$-modules.

\subsection{A complex structure on quantum irreducible flag varieties}
\label{ssect.q.flag.vars}

In \cite{HK04}, \cite{HK06}, \cite{HK07},   Kolb and Heckenberger construct an integrable almost complex structure on the quantum group analogues of irreducible flag varieties. 
 Stefan Kolb kindly provided us with some notes about that work and has allowed us to use them in writing this section. We thank him.
 Of course, any inaccuracies below are due to the present authors.

%The other source for the material in this section is Stokman's paper \cite{St03}.

\subsubsection{Generalized flag varieties}
\label{sect.classical}

Let $\fg$ be a finite dimensional, complex, simple  Lie algebra, $\fp\subset \fg$ 
a parabolic subalgebra, $\fp=\fl \oplus \fn$ its Levi decomposition, and $\fp^{\rm op}:=
\fl \oplus \fn^-$ the opposite parabolic. 
 
Let $G$ be the connected, simply connected, affine algebraic group with Lie algebra $\fg$,
$P\subset G$ the connected subgroup with Lie algebra $\fp$, and $L$ the Levi factor of $P$. 
The quotients $G/L$ and $G/P$ are  affine and projective algebraic varieties, respectively.

\begin{quote}
We make the additional assumption that $\fg/\fp$ is irreducible as a 
$\fp$-module under the adjoint action. In this situation 
$G/P$ is called an  {\it irreducible flag manifold}.\footnote{For a classification of such pairs $(\fg,\fp)$ see \cite[Sect. 3.1]{BE89}.} Irreducible flag manifolds coincide with compact Hermitian
symmetric spaces \cite[Sect. X.6.3]{Helg}. 
\end{quote}

\noindent

% Equivalently, the spaces $G/P$ are the irreducible compact Hermitian symmetric spaces. A list of these manifolds can be found in \cite[Section X.6.3]{b-Helga78} or \cite[p.~27]{b-BastonEastwood}.

Suppose that $\mu$ is a dominant integral weight, $V(\mu)$ is the finite dimensional 
simple $U(\fg)$-module with highest weight $\mu$, and $v_\mu$ is a highest weight 
vector  in $V(\mu)$.

For each parabolic subgroup $P \subset G$ there is a weight $\l$ such that $G/P$ is isomorphic to the $G$-orbit in $\PP\big(V(\l)\big)$ of $\overline{v_\l}$, the image of $v_\l$.\footnote{The weight $\l$ is not unique but there is a ``smallest'' such $\l$, the one that is a sum of fundamental weights with coefficients 0 or 1. From now on it is that $\l$ which we will use.}
%Thus $V(\l)^*$ is linear forms on $\PP\big(V(\l)\big)$.   
The vector space 
 $$
 S[G/P]:=\bigoplus_{n=0}^\infty V(n \l)^*
 $$
endowed with the Cartan multiplication and $\Z$-grading given by $\deg V(\l)^*=+1$ is a homogeneous coordinate ring of $G/P$. It coincides with the subalgebra of $\cO(G)$ generated 
by the matrix coefficients $\{c^\l_{f,v_\l} \; | \;  f \in V(\l)^*\}$ where
$$
c^\mu_{f,v}(u):=f(uv), \qquad f \in V(\mu)^*, \; v \in V(\mu), \; u \in U(\fg).
$$
There are invertible $\cO_{G/P}$-modules $\cL_{n\l}$ such that $H^0(G/P,\cL_{n\l}) \cong V(n \l)^*$.

\subsubsection{Quantum group analogues}

Let $q \in \RR-\{0\}$.

Let $U_q(\fg)$ be the quantized enveloping algebra of $\fg$ defined at \cite[6.1.2]{KS}. The quantum group $G_q$ is defined implicitly by declaring that its coordinate ring, $\cO(G_q)$, is the linear span of the matrix coefficients $c^\mu_{f,v} \in U_q(\fg)^*$, 
$$
c^\mu_{f,v}(u):=f(uv), \qquad f \in V(\mu)^*, \; v \in V(\mu), \; u \in U_q(\fg),
$$
as $\mu$ ranges over the dominant integral weights and $V(\mu)$ is the simple $U_q(\fg)$-module with highest weight $\mu$.

Following Soibelman \cite{S92}, and subsequent authors, let
 $S[G_q/P]$  be the subalgebra  of  $\cO(G_q)$ generated by $\{c^\l_{f,v_\l} \; | \; f \in V(\l)^*\}$.
 There is an isomorphism
 $$
 S[G_q/P] \cong \bigoplus_{n=0}^\infty V(n\l)^*
 $$
 of vector spaces that becomes an algebra isomorphism when the right-hand side is given a quantized analogue of the Cartan product. 
We view $S[G_q/P]$ as a homogenous coordinate ring of a non-commutative analogue of
$G/P$.\footnote{We are writing $S[G_q/P]$ for the algebra denoted $S_q[G/P]$ by Heckenberger and Kolb \cite[Sect. 3.1.1]{HK06} and, later, will write $S[G_q/P^{\rm op}]$ for the algebra they denote by $S_q[G/P^{\rm op}]$.}
In this setting, $S[G_q/P]$ is the algebra denoted $R$ in section \ref{ssect.defn.X}, and the non-commutative analogue of $G \cdot \overline{v_\l}$ plays the role of $X$. 

So far we are only thinking of $G_q/P$ as a non-commutative complex projective variety, $\Projnc(S[G_q/P])$ where $S[G_q/P]$ is thought of as holomorphic polynomial functions on the ``punctured cone'' over $G_q/P$, i.e., the non-commutative space ${\bf X}-\{0\}$ in \S\ref{ssect.defn.X}. 
To give $G_q/P$, or, more precisely, ${\bf X}-\{0\}$, an underlying real structure we need a $*$-algebra that plays the role played by ${\bf R}$ in section \ref{sect.real.R}. 

At a minimum, we need an algebra to play the role of anti-holomorphic polynomial 
functions on ${\bf X}-\{0\}$. The $G$-orbit in $\PP(V(\l)^*)$  
of a lowest weight vector $f_{-\l} \in V(\l)^* \cong V(-w_0\l)$ is isomorphic to 
$G/P^{\rm op}$ so, following the classical case (e.g., \cite[\S\S 2-3]{HK04}),
we define the ring of anti-holomorphic polynomial functions on $G_q/P$ to be the subalgebra $S[G_q/P^{\rm op}] \subset \cO(G_q)$ generated by $\{c^{-w_0\l}_{v,f_{-\l}} \; | \; v \in V(-w_0\l)^* \cong V(\l)\}$. 
It is isomorphic to 
 $$
 \bigoplus_{n=0}^\infty V(n \l)
 $$
endowed with the quantized Cartan multiplication. We consider $S[G_q/P^{\rm op}]$ as a $\Z$-graded algebra concentrated in degree $\le 0$ with $\deg V(\l)=-1$. 
  
 \subsubsection{The underlying real structure on $G_q/P$}
 
 First, we need an underlying real structure on $G_q$ by viewing it as the complexification of a non-commutative real algebraic variety. Recall, $G$ is the complexification of its compact connected real form, $U \subset G$ and, if $K=U \cap P$, the natural inclusion $U/K \to G/P$ is a diffeomorphism.

The Hopf $*$-algebra structure on $U_q(\fg)$ induces a Hopf $*$-algebra structure on 
$\cO(G_q)$. The explicit form of $*$ appears at \cite[(3.5)]{SD99}.
We write $\C[U_q]$ for $\cO(G_q)$ viewed as a $*$-algebra; as a $*$-algebra 
$\C[U_q]$ behaves as $\C$-valued polynomial functions on the quantum group analogue $U_q$ of the compact real form $U$. 
The involution $*$ is an anti-isomorphism between $S[G_q/P]$ and $S[G_q/P^{\rm op}]$ (see \cite[Defn.2.2]{St03}).

Let ${\bf R}$ be the $\Z$-graded vector space $S[G_q/P] \otimes S[G_q/P^{\rm op}]$ with the
following algebra structure: let $v_1,\ldots,v_N$ and $f_1,\ldots, f_N$ be dual bases for $V(\l)$ and $V(\l)^*$ consisting of weight vectors; declare that the subalgebra of ${\bf R}$ generated by $f_1,\ldots,f_N$ (resp., $v_1,\ldots,v_N$) is isomorphic to $S[G_q/P]$ (resp., $S[G_q/P^{\rm op}]$), and declare that $(1 \otimes v_i) (f_i \otimes 1)$ is the linear combination of the elements $f_k \otimes v_\ell$ given by \cite[(28)]{HK06}. 

The degree-zero component, ${\bf R}_0$, can also be described in the following alternative ways:
\begin{enumerate}
  \item 
${\bf R}_0=\AA^q_\l$, the unital $*$-subalgebra of $\C[U_q]$ generated by  the products 
 $ \{ f\cdot g^*  \; | \; f,g \in V(\l)^*\}$ (\cite[Defn. 2.3(b)]{St03});
  \item 
${\bf R}_0=\C[U_q/L\cap U_q]=\C[G_q/L]$, the subalgebra $\{\phi \in \C[U_q] \; | \; u.\phi=\ve(u)\phi \;\; \forall u \in U_q(\fl)\}$  of $\C[U_q]$ (\cite[(4.4)]{SD99}).
\end{enumerate}
Equality of the algebras in (1) and (2) is proved at \cite[Thm. 4.10]{SD99} and \cite[Thm. 2.5]{St03} by passing to C$^*$-algebra completions, and at \cite[Prop. 2]{HK04} without using C$^*$-algebras.

Following Heckenberger and Kolb \cite[\S 3.1.2]{HK06} there is a central element and an associated quotient algebra,
$$
c:=\sum_{i=1}^N v_i \otimes f_i \in {\bf R}_0
\qquad \hbox{and} \qquad 
\overline{\bf R}:=\frac{{\bf R}}{(c-1)}.
$$ 
Because $\deg(c)=0$, $\overline{\bf R}$ inherits the $\Z$-grading on ${\bf R}$. Because $c$ is self-adjoint it also inherits the $*$-algebra structure on ${\bf R}$. The algebra $\overline{\bf R}$ is denoted $S_q[G/P]^{{}^{c=1}}_\C$ by Heckenberger and Kolb in \cite[\S 3]{HK06}.

\subsubsection{The de Rham complex for $\C[G_q/L]$}
\label{sect.de.Rham}

By \cite[Thm.2, p.483]{HK04}, there are exactly two non-isomorphic finite dimensional irreducible 
left covariant first order differential calculi, % (FODC), 
$\C[G_q/L] \to \Gamma_+$ and $\C[G_q/L] \to\Gamma_-$,  
in the sense of Woronowicz \cite[Sect. 12.1.2]{KS}.  
In \cite[Sect. 7.1]{HK07}, Heckenberger and Kolb relabel these  $(\Gamma_\pd,\pd)$
and $(\G_{\pdol},\pdol)$. To conform with our notation we label them 
$(\Omega^{1,0}\C[G_q/L],\pd)$ and $(\Omega^{0,1}\C[G_q/L],\pdol)$. Heckenberger and Kolb now define
$$
(\Omega^1\C[G_q/L],\extd):=(\Omega^{1,0}\C[G_q/L] \oplus \Omega^{0,1}\C[G_q/L], \pd \oplus \pdol).
$$
Thus, $\Omega^1\C[G_q/L]$ and the  derivation $\extd:\C[G_q/L] \to \Omega^1\C[G_q/L]$ provide a 
non-commutative analogue of the module of K\"ahler differentials for the affine variety $G/L$.

Associated to each first order calculus $(\Omega^1\C[G_q/L],\extd)$, $(\Omega^{1,0}\C[G_q/L],\pd)$,
and $(\Omega^{0,1}\C[G_q/L],\pdol)$, is a unique universal higher order differential calculus (see, e.g.,  \cite[12.2.3]{KS}) that we denote by $(\Omega^\hdot\C[G_q/L],\extd)$, $(\Omega^{\hdot,0}\C[G_q/L],\pd)$,
and $(\Omega^{0,\hdot}\C[G_q/L],\pdol)$, respectively.\footnote{These are denoted 
%$(\G_{d,u}^\wedge,\extd)$,  $(\G_{\pd,u}^\wedge,\pd)$, and  
%$(\G_{\pdol,u}^\wedge,\pdol)$ in \cite[\S 3.3]{HK06} and  
$(\G_d,\extd)$,  $(\G_\pd,\pd)$, and $(\G_{\pdol},\pdol)$, in the sentence before \cite[Thm.7.1]{HK07}.} 

For example, $(\Omega^\hdot\C[G_q/L],\extd)$  is the differential graded algebra that is
the quotient of the tensor algebra $T_{\C[G_q/L]}(\Omega^1)$ by the ideal generated by 
\begin{equation}
\label{univ.diffl.calc}
\Bigg\{\sum_i\extd a_i \otimes \extd b_i \; \Bigg \vert \; \sum_ia_i \extd b_i=0\Bigg\}
\end{equation}
 with differential defined by 
$\extd(a_0 \extd a_1 \wedge \ldots \wedge \extd a_n)=\extd a_0 \wedge \extd a_1 \wedge \ldots \wedge \extd a_n$.\footnote{Here, as throughout this paper, $\wedge$ is used to denote the multiplication and is not assumed to be skew-commutative.} 
It follows from the definition of $\Omega^\hdot\C[G_q/L]$ and $\extd:\Omega^\hdot\C[G_q/L]\to \Omega^{\hdot+1}\C[G_q/L]$
that $\extd^2=0$.

Heckenberger and Kolb call $(\Omega^\hdot\C[G_q/L],\extd)$ the de Rham complex for $\C_q[G/L]$.
By construction, $\Omega^0\C[G_q/L]=\C_q[G/L]$ and $\extd:\C_q[G/L] \to \Omega^1\C[G_q/L]$ is, as above, $\pd \oplus \pdol$.  

The main results in \cite{HK06}, Propositions 3.6, 3.7, and 3.11, show that the homogeneous components of $(\Omega^\hdot\C[G_q/L],\extd)$, $(\Omega^{\hdot\C[G_q/L],0},\pd)$,
and $(\Omega^{0,\hdot}\C[G_q/L],\pdol)$, have the same dimensions as their classical counterparts. 
In particular, $\Omega^{2\dim(\fg/\fp)}\C[G_q/L]$ is a free left and right $\C[G_q/L]$-module of rank one, a generator being interpreted as a volume form on $G_q/P$. 

\subsubsection{$(p,q)$-forms on $G_q/P$}
In \cite{HK07}, it is shown that $(\Omega^\hdot\C[G_q/P],\extd)$ may be constructed in terms of the 
Bernstein-Gelfand-Gelfand resolution and that that construction  \cite[above Theorem 7.14]{HK07} 
gives a decomposition of $(\Omega^\hdot\C[G_q/P],\extd)$ as a direct sum of covariant $\C_q[G/L]$-bimodules 
\begin{align}\label{sum}
  \Omega^n\C[G_q/P]=\bigoplus_{p+q=n}\Omega^{p,q}\C[G_q/P]
\end{align}
where $\Omega^{1,0}\C[G_q/P]$ and $\Omega^{0,1}\C[G_q/P]$ are as in section \ref{sect.de.Rham} and 
 $\extd:\Omega^n\C[G_q/P]\rightarrow \Omega^{n+1}\C[G_q/P]$ satisfies $\extd=\pd+\pdol$ with
$\pd(\Omega^{p,q}\C[G_q/P])\subset\Omega^{p+1,q}\C[G_q/P]$ and $\pdol(\Omega^{p,q}\C[G_q/P])\subset\Omega^{p,q+1}\C[G_q/P]$. The differentials $\pd$ and $\pdol$ coincide with those defined in \cite[Proposition 3.8, Remark 3.10]{HK06}. Proceeding as in \cite[Lemmata 7.6, 7.17]{HK07}, one shows that the multiplication maps
\begin{align*}
  \wedge: \Omega^{1,0}\C[G_q/P]: \Omega^{p,q}\C[G_q/P]\rightarrow \Omega^{p+1,q}\C[G_q/P]\\
  \wedge: \Omega^{0,1}\C[G_q/P]: \Omega^{p,q}\C[G_q/P]\rightarrow \Omega^{p,q+1}\C[G_q/P]
\end{align*}
are surjective, and hence
\begin{align}\label{mult-compatible}
  \Omega^{p,q}\C[G_q/P]\wedge \Omega^{p',q'}\C[G_q/P]=\Omega^{p+p',q+q'}\C[G_q/P]\ .
\end{align}

\subsubsection{Definition of $J$}

Let $J:\Omega^1\C[G_q/P] \to \Omega^1\C[G_q/P]$ be the unique $\C$-linear map such that $J^2=-1$ and 
$\Omega^{1,0}$ is the $+1$-eigenspace and $\Omega^{0,1}$ is the $-1$-eigenspace.
It is clear that $J$ is a $\C_q[G/L]$-bimodule automorphism of both $\Omega^{1,0}\C[G_q/P]$ and 
$\Omega^{0,1}\C[G_q/P]$ hence a $\C_q[G/L]$-bimodule automorphism of $\Omega^{1}\C[G_q/P]$.

Using the direct sum decomposition (\ref{sum}) we can extend $J$ to a bimodule isomorphism $J:\Omega^\hdot\C[G_q/P]\rightarrow \Omega^\hdot\C[G_q/P]$ by
\begin{align*}
  J(\omega)=(p-q)\mathrm{i} \omega \qquad \mbox{for all $\omega\in \Omega^{p,q}\C[G_q/P]$}.
\end{align*}
By (\ref{mult-compatible}), $J$ is a derivation: if $\omega\in \Omega^{p,q}\C[G_q/P]$ and $\tau\in\Omega^{p',q'}\C[G_q/P]$, then 
\begin{align*}
  J(\omega\wedge\tau)&=(p+p'-(q+q'))\mathrm{i}\,\omega\wedge\tau\\
                     &=(p-q)\mathrm{i}\,\omega\wedge\tau + (p'-q')\mathrm{i}\,\omega\wedge\tau\\
                     &=J(\omega)\wedge \tau + \omega\wedge J(\tau).
\end{align*}

\end{document}